\documentclass[twocolumn]{autart}    % Enable this line and disable the
\usepackage{graphics} % for pdf, bitmapped graphics files
\usepackage{epsfig} % for postscript graphics files
\usepackage{mathptmx} % assumes new font selection scheme installed
\usepackage{times} % assumes new font selection scheme installed
\usepackage{amsmath} % assumes amsmath package installed
\usepackage{amssymb}  % assumes amsmath package installed
\usepackage{color}
\usepackage{graphicx,colordvi}

\newtheorem{lemma}{ Lemma}
\newtheorem{theorem}{ Theorem}

\newtheorem{remark}{ Remark}
\newtheorem{example}{ Example }
\newtheorem{corollary}{ Corollary}
\newtheorem{proposition}{ Proposition}          % Include this line if your
\begin{document}

\begin{frontmatter}
%\runtitle{Insert a suggested running title}  % Running title for regular
                                              % papers but only if the title
                                              % is over 5 words. Running title
                                              % is not shown in output.

\title{Mean Field Production Output Control with Sticky Prices: Nash and Social Solutions\thanksref{footnoteinfo}} % Title, preferably not more
                                                % than 10 words.

\thanks[footnoteinfo]{This paper was not presented at any IFAC
meeting. Corresponding author Bingchang Wang.}

\author[BW]{Bingchang Wang}\ead{bcwang@sdu.edu.cn},\    % Add the
\author[MH]{Minyi Huang}\ead{mhuang@math.carleton.ca}               % e-mail address
%\author[Baiae]{Publius Maro Vergilius}\ead{vergilius@culture.ir}  % (ead) as shown

\address[BW]{School of Control Science and Engineering, Shandong
        University, Jinan 250061, P. R. China}  % Please supply
\address[MH]{School of Mathematics and Statistics, Carleton University, Ottawa, ON K1S 5B6, Canada}             % full addresses

                             % keyword list or with the
                                          % help of the Automatica
                                          % keyword wizard

\begin{abstract}                          % Abstract of not more than 200 words.
This paper presents an application of mean field control to  dynamic production optimization. Both noncooperative and cooperative solutions are considered. We first introduce a market of a large number of agents (firms) with sticky prices and adjustment costs. %The cost functional of each agent has indefinite state weights.
By solving auxiliary limiting optimal control problems subject to consistent mean field approximations, two sets of decentralized strategies are obtained and further shown to asymptotically attain Nash equilibria and social optima, respectively. The performance estimate of the social optimum strategies exploits a passivity property of the underlying model. A numerical example is given to compare market prices, firms' outputs and costs under two two solution frameworks.

\end{abstract}

\end{frontmatter}

\section{Introduction}

Mean field game theory is effective to
design decentralized strategies in a system of many players which
are individually negligible but collectively affect a particular
player (see e.g., \cite{HCM03}, \cite{HCM07}, \cite{HMC06}, \cite{LL07}).
 By identifying a consistency relationship between the individual's best response  and the mass (population macroscopic) behavior, one may
 obtain a fixed-point equation to specify the mean field.  This procedure leads to a set of  decentralized  strategies as an $\epsilon$-Nash equilibrium for the actual model with a  large but finite population.
By now, mean field games have been intensively studied in the LQG (linear-quadratic-Gaussian) framework \cite{HCM03}, \cite{HCM07}, %\cite{KC13},
\cite{ELN13}, \cite{TZB14}, \cite{FMM16}; there is also a
large body of works on nonlinear models \cite{HMC06}, \cite{LZ08A}, %\cite{CD13},
 \cite{BFY13}, \cite{GS13}. %\cite{NC13}.  BFY13, GS13}.
For further
literature, readers are referred to \cite{huang2010large}, %\cite{NH12},
 \cite{wang2012major}, \cite{WZ14} for mean field
models with a major player, \cite{weintraub2008markov} for
 oblivious equilibria proposed for large-scale Markov decision processes
  of industry dynamics,
\cite{WZ12}
for mean field games with Markov jump parameters. %\cite{YMMS12} for mean field oscillator games, \cite{AJW13} for a solution notion called stationary equilibrium,
% \cite{CS14} for application to Bertrand and Cournot equilibrium models.
For a survey
on mean field game theory, see \cite{BFY13}, \cite{GS13}, and \cite{CHM17}. Besides noncooperative games, social optima in mean field control have been investigated in some literature \cite{HCM12}, \cite{WZ17}.
Mean field games and control
have found wide applications, including smart grids  \cite{MCH13}, \cite{CBM15}, \cite{KM16}, finance, economics \cite{weintraub2008markov}, \cite{GLL11}, \cite{CS14}, \cite{HN16}, operation research \cite{LW11}, \cite{AJW13}, \cite{LM14}, \cite{SMJ15}, and social sciences \cite{LT15}, \cite{HM16}, \cite{BTN16}, etc.

%Thus, at
%each point in time oligopolists face a constant price that they know will
%decline, but not instantaneously, if their joint output exceeds the level of demand
%at that price.
%The evolution of price through time is a function of the difference
%between the current price and the price on the demand function for each level
%of.output.

This paper aims to present an application of mean field control to production output adjustment in a large market with many firms and sticky prices. Under the stickiness assumption, the price of the underlying product does not adjust instantaneously according to its demand function, %at the given level of output,
but evolves slowly and smoothly. Dynamic game models for duopolistic competition with sticky prices were initially proposed by Simaan and Takayama \cite{simaan1978game}, and then extended to investigate asymptotically stable steady-state equilibrium prices
in \cite{fershtman1987dynamic}.
In \cite{CR04}, \cite{wiszniewska2014dynamic}, the authors considered open and closed-loop Nash equilibria for dynamic oligopoly with $N$ firms and compared prices' behavior in and outside the steady-state levels, respectively. \Blue{Adjustment costs in production models have been addressed in the economic
 literature (see e.g. \cite{R02}) and {they have been taken into account in the study of dynamic oligopoly \cite{driskill1989dynamic}, %\cite{Don92},
  \cite{Sc1997}, \cite{JV04}.}
The work \cite{driskill1989dynamic} introduces a duopoly where each firm has output level %$q_i$
 subject to control according to a first-order integrator
 dynamics.
 However, when the number of firms is large (e.g. in a perfectly competitive market) and the adjustment cost is considered, the computational complexity of output adjustment is high.}
In the mean field control framework, one can effectively address the complexity issue.

  Within our model, a large number of producers supply a certain product with sticky prices, and the output adjustment incurs a cost. The cost function of a firm is based on product cost, price, and adjustment cost. In \cite{WH2015}, we combined the price and firm's output as a 2-dimensional system. Thus, the cost function has indefinite state weights, which differs from many existing LQG models of mean field games in the literature
  \cite{HCM07}, \cite{WZ12}. In this paper, the price in the mean field limit model is taken as an exogenous signal without the need of state space augmentation. This contributes to deriving a simple condition that ensures the solvability of the resulting equation system.

The Nash equilibrium and the social optimum are two fundamental solution notions to competitive markets with many firms, where the former applies to the noncooperative model,
and the latter is for
 the cooperative model.
In this paper, we design Nash and social optimum strategies for the production output control model based on the mean field control methodology, respectively, and further compare two solutions numerically.
%We first present the Nash solution within the noncooperative game formulation, and then give the social solution, where agents cooperate to optimize the social %cost.
The Nash solution of our model starts by solving a limiting optimal control problem and next applies the consistency requirement for the mean field approximation.
 We then obtain a set of decentralized strategies and show that the set of strategies is an $\varepsilon$-Nash equilibrium. For the social optimum solution, we first provide an auxiliary optimal control problem by a person-by-person optimality approach, and then design
a set of decentralized strategies by solving the limiting auxiliary problem subject to consistent mean field approximations.  The set of strategies is shown to be asymptotically socially optimal by exploiting a passivity property of the underlying model.

 An illustrative numerical example is given to compare market prices, firms¡¯ outputs and optimal costs under the
game and social optimum frameworks. It is numerically shown that the social optimum has a
 lower average output level than that in the noncooperative case.
This is similar to the behavior in a duopoly model \cite{V93}
where cooperation of the two players results in
a lower total output than in the Cournot equilibrium.

The paper is organized as follows. Section II introduces the game and social optimum problems with $N$ players. In Section III, we first design a set of decentralized strategies by the mean field control methodology and then show its asymptotic Nash equilibrium property. In Section IV, we construct a set of decentralized strategies, which is shown to be asymptotically socially optimal. In Section V, a comparison of two solutions is demonstrated by a numerical example.
%we construct the equation system for the solution of the mean field game, and analyze its solution. In Section V, we show that the closed-loop system is stable %and the proposed set of strategies is an $\varepsilon$-Nash equilibrium.
Section VI concludes the paper.

Notation: $\|\cdot\|$
denotes the Euclidean vector norm or matrix spectral norm. For a matrix $M$, $|M|$ denotes the determinant of $M$.
$C([0,\infty),\mathbb{R}^n)$ denotes the class of $n$-dimensional
continuous functions on $[0,\infty)$; $C_b([0,\infty),\mathbb{R}^n)$ is the class of bounded and continuous functions; $C_{\rho}([0,\infty), \mathbb{R}^n)=\{f|f\in C([0,\infty), \mathbb{R}^n), \sup_{t\geq 0}|f(t)|e^{-\rho^{\prime}t}<\infty \ \hbox{for\ some }\
 \rho^{\prime}\in [0,\rho) \}.$
 For a family of $\mathbb{R}^n$-values random variables $\{x(\lambda),
 \lambda\geq 0\}$,
$\sigma(x(\lambda), \lambda\leq t)$ is the $\sigma$-algebra  generated by the collection of random variables; $\|x\|_{\rho}=\big[E\int_0^{\infty} e^{-\rho t} |x(t)|^2d t\big]^{1/2}$; $\|x\|_{\infty}=\sup_{t\geq 0} |x(t)|$.
For two sequences $\{a_n, n = 0, 1, \cdots\}$ and $\{b_n, n = 0, 1,  \cdots\}$, $a_n = O(b_n)$ denotes
$\limsup_{n\to\infty}|{a_n}/{b_n}|\leq C$, and $a_n = o(b_n)$ denotes
$\limsup_{n\to\infty}|{a_n}/{b_n}|=0$.
 For convenience of presentation, we use $C, C_1, C_2,\cdots$ to
denote generic positive constants, which may vary from place to place.

\section{Problem Description}

\subsection{Dynamic oligopoly with sticky prices}

Dynamic game models for
oligopolistic competition with sticky prices were initially
proposed by Simaan and Takayama \cite{simaan1978game}, and then further investigated
in \cite{fershtman1987dynamic}, \cite{CR04}, \cite{wiszniewska2014dynamic}. %The models with two firms are considered in \cite{simaan1978game} and \cite{fershtman1987dynamic}.
%The authors in \cite{CR04, wiszniewska2014dynamic}
%investigated open and closed-loop Nash equilibria for dynamic
%oligopoly with $N$ firms and compared prices¡¯ behavior in
%and outside the steady-state levels, respectively.
According to the model in \cite{CR04}, \cite{wiszniewska2014dynamic}, the sticky price evolves by
$$\frac{dp}{dt}=\alpha(\beta-\delta\sum_{j=1}^Nq_j-p),\quad p
(0) \hbox{ given}, $$
where $q_j$ is the output of firm $j$, $j=1,\cdots,N$, and has the role of control.
The payoff function of firm $i$ is described by
$$K_i(q_1,\cdots,q_N)=\int_0^{\infty}e^{-\rho t}(pq_i-cq_i-\frac{1}{2}q_i^2)dt.$$
The constants $\alpha,\beta,\delta$ and $c$ are positive, and $c$ is the cost of unit output.

%On the other hand, adjustment
%costs are important in some scenarios \cite{R02} and have been taken into account in the study of dynamic oligopoly \cite{driskill1989dynamic}, \cite{JV04}.  %as evidenced by empirical studies at the micro level.
%The work \cite{driskill1989dynamic} developed a duopoly model in which firms incur costs according to how fast they adjust their  output levels.
%Specifically, firm $i$ is to solve the
%problem
%
%$$\max_{u_i} K_i=\int_0^{\infty} e^{-\rho t}[P(q_1,q_2)q_i-C_i(q_i,u_i)]dt,$$
%subject to
%$$\dot{q}_i=u_i,\ i=1,2,$$
%where $q_i$ denotes the output of firm $i$, $u_i$ is its changing rate, and
%$$P=a-(q_1+q_2),$$
%$$C_i=cq_i+\frac{b}{2}q_i^2+\frac{A}{2}u_i^2.$$

\subsection{Output adjustment in a mean field framework}

The paper considers a large market of many firms.
Based on the formulation of sticky prices in \cite{fershtman1987dynamic}, \cite{wiszniewska2014dynamic}, we assume that the price evolves by
\begin{align}\label{eq1}
\frac{d{p}(t)}{dt}=&\hspace*{0.05cm}\alpha[\beta-p(t)-q^{(N)}(t)]\cr
=&-\alpha p(t)-\alpha q^{(N)}(t)+\alpha\beta,
\end{align}
where $\alpha>0$ denotes the speed of adjustment to the level on the demand function,
%\beta-p(t)-q^{(N)}(t) is the price on the demand function for the given level of  firms' production.
and $q^{(N)}(t)=\frac{1}{N}\sum_{i=1}^Nq_i(t)$ is the average of firms' outputs.
The output of each firm is described by the stochastic differential equation (SDE)
\begin{equation}\label{eq2}
  dq_i(t)=-\mu q_i(t)dt+b_iu_i(t)dt+\sigma dw_i(t),
\end{equation}
where $\{w_i(t), i=1,\cdots,N\}$ are independent standard Brownian motions, which are also independent of initial outputs of all firms $\{q_i(0), i=1,\cdots,N\}$.
The constants $\alpha,\beta$, $\mu$ and $b_i$ are positive.

Adjustment costs in production models have been addressed in the economic
 literature (see e.g. \cite{R02}) and {they have been taken into account in the study of dynamic oligopoly \cite{driskill1989dynamic}, \cite{Don92}, \cite{Sc1997}, \cite{JV04}.}
The work \cite{driskill1989dynamic} introduces a duopoly where each firm has output level $q_i$
 subject to control $u_i$ according to a first order integrator
 dynamics. In the resulting differential game, the instantaneous payoff of each firm
 is determined from its net profit minus  quadratic penalty terms of
 $q_i$ and $u_i$

\begin{remark}
{As in \cite{fershtman1987dynamic}, $\beta-q^{(N)}$ is the price on the demand function for the given level of firms' outputs.}
In the static case, the inverse demand function has a linear version
$p = \beta - \delta q^{(N)}$; here for simplicity we set $ \delta$ as 1. The scaling factor $1/N$ for $q^{(N)}$ is standard in modelling
and analysis of large markets, and some closely related price modelling in a large dynamic market can be
found in \cite{L84}, \cite{UY89}. $\mu$ is used to indicate friction in adjusting the output, and
$w_i$ is random shocks in output.
\end{remark}

The cost function of each firm is given by
\begin{equation}\label{eq3}
J_i(u)=E\int_0^{\infty} e^{-\rho t}L(p,q_i,u_i)dt,
\end{equation}
where $$L=-p(t)q_i(t)+cq_i(t)+ru_i^2(t),$$
$$u=(u_1,\cdots,u_{i},\cdots,u_N),$$
$r>0$ and $0<c<\beta$.
  Here, $c$ denotes the production cost,
  and $ru_i^2(t)$ denotes the adjustment cost.  The minimization of $J_i(u)$
is equivalent to maximizing the payoff
\begin{equation*}
K_i(u)=E\int_0^{\infty} e^{-\rho t}[q_i(t)(p(t)-c)-ru_i^2(t)]dt.
\end{equation*}
{We only consider the case $\beta> c$ to make the subsequent optimization problems be of practical interest. Otherwise, given a positive $q^{(N)}$,
the production cost already exceeds the price, and the optimization problem is not too meaningful.  }

The social cost is defined as
\begin{equation}\label{eq4}
%J^{(N)}(u)=\sum_{j=1}^N J_j(u).
J^{(N)}_{\rm soc}(u)=\frac{1}{N}\sum_{j=1}^N J_j(u).
\end{equation}
%\begin{remark}

%The assumption $\beta> c$ ensures that the output in the steady state is positive.
%\end{remark}

Based on costs (3) and (4), one may formulate a standard LQG game and an optimal control problem, respectively. A limitation of this approach
is that the control strategy will be centralized.
Our goal is to look for decentralized strategies for the corresponding optimization criterion.

The basic objective of this paper is to seek Nash solutions and social solutions to
mean field production output control with sticky prices. Specifically, we study the following two problems:

\emph{Problem I:} Find $\varepsilon$-Nash equilibrium strategies for agents to minimize the individual cost $J_i$
over the set of decentralized strategies
$${\mathcal U}_{d,i}=\Big\{u_i: u_i(t)\ \hbox{is adapted to} \  {\mathcal F}_t^i, E\int_0^{\infty} e^{-\rho t} u_i^2(t)dt<\infty%\\hbox{such that}\ E[p^2(t)+q_i^2(t)]=o\big(e^{\rho t}\big)
\Big\},$$ where $ {\mathcal F}^i_t=\sigma\{q_i(0),w_i(s), s\leq t\}$, $t\geq 0$, $i=1,\cdots,N$.

\emph{Problem II:} Find asymptotic social optimum strategies for agents to minimize $J^{(N)}_{\rm soc}$
over the set of decentralized strategies
${\mathcal U}_{d,i}$, $i=1,\cdots,N$.

\Blue{For a large market, a natural way of modeling the sequence of parameters $b_1,\cdots,b_N$ is to view them as being sampled from a space such that this sequence exhibits certain statistical properties when $N\to\infty$. Define the associated empirical distribution function $F_N(\theta)=\frac{1}{N}\sum_{i=1}^NI_{[\theta_i\leq \theta]}$, where $I_{[\theta_i\leq \theta]}=1$ if $\theta_i\leq \theta$ and $I_{[\theta_i\leq \theta]}=0$ otherwise.}

%The objective for each agent is to minimize $J^{(N)}(u)$ over the decentralized set
%$${\cal U}_{d,i}=\Big\{u_i: u_i(t)\ \hbox{is adapted to} \  {\cal F}_t^i%,\\hbox{such that}\ E[p^2(t)+q_i^2(t)]=o\big(e^{\rho t}\big)
%\Big\},\quad i=1,\cdots,N,$$ where $ {\cal F}^i_t=\sigma\{q_i(0),w_i(s), s\leq t\}, t\geq 0$.

%

%Different from \cite{CR04, wiszniewska2014dynamic}, the problem (\ref{eq1})-(\ref{eq3}) involves changing rate of output and adjustment
%costs. In this case, the output levels are state variables which do not jump but evolve continuously. Adjustment
%costs are important in some scenarios \cite{R02} and have been taken into account in the study of dynamic oligopoly \cite{driskill1989dynamic, JV04}.  %as evidenced by empirical studies at the micro level.
%The work \cite{driskill1989dynamic} developed a duopoly model in which firms incur costs according to how fast they adjust their level of output.

We introduce the assumptions.

\textbf{A1)} The initial price $p(0)=p_0>0$ is a constant. The initial outputs of all firms $\{q_i(0), i=1,\cdots,N\}$ are independent. $Eq_i(0)=q_0>0$ for all $i=1,\cdots,N$; there exists $c_0<\infty$ independent of $N$
such that $\max_{i=1,\cdots,N}E|q_i(0)|^2\leq c_0$.

\textbf{A2)} There exists a distribution function $F(\cdot)$ such that the empirical distribution $F_N$ converges weakly to $F$, where $F_N(b)=\frac{1}{N}\sum_{i=1}^NI_{[b_i\leq b]}$. Furthermore, each
$b_i>0$ and $\int_{\mathbb{R}}\theta^2dF(\theta)>0$.

\textbf{A3)} For all $N$, $\{b_i,i=1,\cdots,N\}$ is contained in a fixed compact set $\Theta$,
 and $\int_{\Theta}\theta dF(\theta)=1$.

\section{Nash Solutions to Output Adjustment} %Production
\subsection{Optimal control for the limiting problem}

Assume that $\bar{q}\in C_{\rho/2}([0,\infty), \mathbb{R})$ is given for approximation of $q^{(N)}$.
Replacing $q^{(N)}$ in (\ref{eq1}) by $\bar{q}$, we introduce
\begin{equation}\label{eq5a}
\frac{d{\bar{p}}(t)}{dt}=\alpha[\beta-\bar{p}(t)-\bar{q}(t)],  \quad \bar{p}(0)=p_0.
\end{equation}
 Accordingly, by replacing $p$ in (\ref{eq3}) with $\bar{p}$ we define the cost function:
\begin{equation}\label{eq5}
\bar{J}_i(u_i)=E\int_0^{\infty} e^{-\rho t}[ (c-\bar{p})q_i
+ru_i^2]dt.
\end{equation}
%Denote $ {\cal F}_t^{i}=\sigma\{p(0),q_i(0),w_i(s), s\leq t\}, t\geq 0$.
{The corresponding admissible control set is ${\mathcal U}_{d,i}$.}
%$$\bar{{\mathcal U}}_i=\big\{u_i: u_i\ \hbox{is adapted to}\  {\mathcal F}^{i}_t,\ \hbox{and}\ E[q_i^2]=o\big(e^{\rho t}\big) \big\}.$$

We first take $\bar{p}$ as an exogenous signal and solve the problem in (\ref{eq2}), (\ref{eq5a}) and (\ref{eq5}).
%Denote
%$$A=\left[\begin{array}{cc}
% -\alpha&0\\0&-\mu
%\end{array}\right], \ B_i= \left[\begin{array}{c}
%  0\\b_i
%\end{array}\right],\ \ Q=\left[\begin{array}{cc}
% 0&-\frac{1}{2}\\-\frac{1}{2}&0
%\end{array}\right],  $$
%and $\bar{x}_i(t)=[\bar{p}(t),q_i(t)]^T$.
For a general initial condition $q_i(t)=q_i$ at time $t$, define the value function
$$V_i(t,{q}_i)=\inf_{u_i\in {\mathcal U}_{d,i}}E\left[\int_t^{\infty}e^{-\rho(\tau-t)}L(q_i,u_i)d\tau\Big|q_i(t)=q_i\right].$$
%By the Bellman's Principle of Optimality, we have
%\begin{eqnarray*}
%  v_i(t,\bar{x}_i)&\hspace*{-0.2cm}=\hspace*{-0.2cm}&\inf_{u_i\in {\mathcal U}_i}E\left[\int_t^{t+\Delta t}e^{-\rho(\tau-t)}L(\bar{x}_i,u_i)d\tau\right. \cr
%  &&\left.+e^{-\rho \Delta t}v_i(t+\Delta t,\bar{x}_i(t+\Delta t) )|\bar{x}_i(t)=\bar{x}_i\right].
%  \end{eqnarray*}
%By using Ito's formula, we get the HJB equation:
We introduce the HJB equation:
\begin{align}\label{eq6n}
 \rho V_i
 =\inf_{u_i\in \mathbb{R}}&\left\{ \frac{\partial V_i}{\partial t}
 +\frac{\partial V_i}{\partial q_i}\left(-\mu q_i+b_iu_i\right) \right.\cr
& \ \left. +\frac{\sigma}{2}\frac{\partial^2 V_i}{\partial q_i^2 } +(c-\bar{p})q_i+ru_i^2\right\},
\end{align}
\Blue{where $V_i\in C_{\rho}([0,\infty),\mathbb{R})$.} Let $V_i(t,q_i)=k_iq_i^2+2s_i(t)q_i+g_i(t)$. Then the optimal control law is
 \begin{equation}\label{eq7n}
   \bar{u}_i=-\frac{1}{2r}b_i\frac{\partial V_i}{\partial q_i}=-\frac{b_i}{r}(k_iq_i+s_i).
  \end{equation}
Substituting the control (\ref{eq7n}) into (\ref{eq6n}),
we obtain
\begin{eqnarray*}
  \rho(k_iq^2_i+2s_i{q}_i+g_i)
&=&(-2\mu k_i-\frac{b_i^2}{r}k_i^2)q_i^2\cr
 &\hspace*{-0.2cm}+\hspace*{-0.2cm}&\hspace*{-0.1cm}2\left[ \frac{ds_i}{dt}+(-\mu-\frac{b_i^2}{r}k_i)s_i+\frac{c-\bar{p}}{2}\right]q_i\cr
&\hspace*{-0.2cm}+\hspace*{-0.2cm}&\hspace*{-0.1cm}\frac{dg_i}{dt}-\frac{b^2_i}{r}s_i^2+\sigma k_i.
\end{eqnarray*}
 This yields
\begin{align}\label{eq6}
\rho k_i=&-2\mu k_i-\frac{b_i^2}{r}k_i^2,\\
\label{eq7a}
  \rho s_i=&\frac{ds_i}{dt}+(-\mu-\frac{b_i^2}{r}k_i) s_i+
 \frac{c-\bar{p}}{2},
\\
 \rho g_i=&\frac{dg_i}{dt}-\frac{b^2_i}{r}s_i^2+\sigma k_i. \label{eq8}
\end{align}
%%------------------------------------------------
%
%%Denote the algebra Riccati equation
%%\begin{equation}\label{eq6}
%%\rho P_i=A^TP_i+P_iA-P_iB_iR^{-1}B_i^TP_i+Q,
%%\end{equation}
%%where $$A=\left[\begin{array}{cc}
%% -\alpha&0\\0&0
%%\end{array}\right], \ B_i= \left[\begin{array}{c}
%%  0\\b_i
%%\end{array}\right],\ R=r,\ Q=\left[\begin{array}{cc}
%% 0&-\frac{1}{2}\\-\frac{1}{2}&0
%%\end{array}\right].  $$
%
%%We have the following lemma.
\begin{lemma}\label{lem1n}
$ k_i=0$ is the unique solution to the algebraic Riccati equation (\ref{eq6}) such that $-\mu-\frac{b_i^2}{r}k_i-\frac{\rho}{2}<0$.
\end{lemma}
\emph{Proof.} By solving (\ref{eq6}), we have $k_i=0$ or $k_i=-\frac{r}{b_i^2}(\rho+2\mu)$. If $k_i=0$,
$-\mu-\frac{b_i^2}{r}k_i-\frac{\rho}{2}=-\mu-\frac{\rho}{2}<0$.
 Otherwise, when $k_i=-\frac{r}{b_i^2}(\rho+2\mu)$, $-\mu-\frac{b_i^2}{r}k_i-\frac{\rho}{2}=\frac{\rho}{2}+\mu>0
$.    \hfill$\Box$

\begin{remark}
  The inequality in Lemma \ref{lem1n} specifies a stability
condition for the
closed-loop system which must be satisfied by the solution of $k_i$.
\end{remark}

\begin{theorem} \label{thm1a}
For the optimal control problem in (\ref{eq2}), (\ref{eq5a}) and (\ref{eq5}), assume that $\bar{q}\in C_{\rho/2}([0,\infty), \mathbb{R})$ is given. Then we have\\
\ \ 1) there exists a unique solution $s_i\in {C}_{\rho/2}([0,\infty), \mathbb{R}) $ to (\ref{eq7a});\\
\ \ 2) the optimal control law is uniquely given by $\bar{u}_i=-\frac{b_i}{r}s_i$;\\
\ \ 3) there exists a unique solution $g_i\in C_{\rho}([0,\infty), \mathbb{R})$ to (\ref{eq8}), and the optimal cost is given by
 $$
V_i(0,q_i(0))=
2s_i(0)q_0+g_i(0).
$$
\end{theorem}
\emph{Proof.} Note that by (\ref{eq5a}), $\bar{q}\in C_{\rho/2}([0,\infty), \mathbb{R})$ implies $\bar{p}\in C_{\rho/2}([0,\infty), \mathbb{R})$.
{We can prove parts 1) and 3) %analyzing the initial values of $s_i$ and $g_i$ and obtaining $s_i(0)$ and $g_i(0)$ are uniquely determined from the fact that $s_i\in
by showing that $s_i(0)$ and $g_i(0)$ are uniquely determined from the fact $s_i\in C_{\rho/2}([0,\infty), \mathbb{R})$ and $g_i\in C_{\rho}([0,\infty), \mathbb{R})$, respectively (see e.g., \cite{huang2010large}, \cite{HCM07})}. {To show part 2) we first obtain a prior integral estimate of $q_i$ (see (\ref{eq12aa})) and then use the completion of squares technique (see e.g., \cite{HCM12}, \cite{WZ17}). By Lemma \ref{lem1a},
 $E\int_0^\infty e^{-\rho t} u_i^2dt <\infty$
 implies  $E\int_0^\infty e^{-\rho t} q_i^2dt <\infty$, which further gives that
 $\bar J_i$ is well defined to be finite since $\bar{p}\in C_{\rho/2}([0,\infty), \mathbb{R})$. By Proposition \ref{prop1},
$\bar J_i(u_i)<\infty$ leads to  $E\int_0^\infty e^{-\rho t} u_i^2dt <\infty$
 which further implies
 \begin{equation}\label{eq12aa}
 E\int_0^\infty e^{-\rho t} q_i^2dt <\infty .
  \end{equation}}
  \hfill$\Box$

\subsection{Control synthesis and analysis}

Following the standard approach in mean field games \cite{HCM03}, \cite{HCM07}, we construct the equation system as follows:
\begin{align}\label{eq9a}
\frac{d{\bar{p}}}{dt}=&\alpha[\beta-\bar{p}-\bar{q}].\\ \label{eq9b}
\rho s=&\frac{ds}{dt}-\mu s+
 \frac{c-\bar{p}}{2}\\ \label{eq9c}
\frac{d\bar{q}_{\theta}}{dt}=&-\mu \bar{q}_{\theta}-\frac{\theta^2}{r}s\\ \label{eq9d}
\bar{q}=&\int_{\mathbb{R}}\bar{q}_{\theta}dF(\theta).
\end{align}

In the above, $\theta$ is a continuum parameter. $\bar{q}_{\theta}$ is regarded as the expectation of the state given the parameter $\theta$ in the individual dynamics. The last equation is due to the consistency requirement for the mean field approximation. $p(0)=p_0,q_{\theta}(0)=q_0$ and $s(0)$ is to be determined.
For further analysis,
we make the following assumption.

 \textbf{A4)} There exists a solution $(s, \bar{q}_{\theta}, \theta\in \mathbb{R})$ to (\ref{eq9a})-(\ref{eq9d}) such that for each $\theta\in \mathbb{R}$, both $s$ and $\bar{q}_{\theta}$ are within $C_{b}([0,\infty), \mathbb{R})$.

Some sufficient conditions for ensuring \textbf{A4)} may be obtained by using the fixed-point methods similar to those in \cite{HCM03}, \cite{HCM07}.

\begin{proposition}
If
$$\begin{aligned}
&\frac{1}{2r\mu(\rho+\mu)}\int_{\Theta}{\theta^2}dF(\theta)<1,
\end{aligned}$$
then \textbf{A4)} holds.$\qquad\qquad\qquad\qquad\qquad\qquad\qquad\qquad\quad \Box$

\end{proposition}

{\emph{Proof.}} By (\ref{eq9a})-(\ref{eq9c}), we have
\begin{align*}
  \bar{p}(t)=&\bar{p}(0)e^{-\alpha t}+\int_0^te^{-\alpha(t-\tau)}[\alpha\beta-\alpha \bar{q}(\tau)]d\tau,\\
   s(t)=&\int_t^{\infty}e^{(\rho+\mu)(t-\tau)}\Big[\frac{\bar{p}(\tau)-c}{2}\Big]d\tau,\\
   \bar{q}_{\theta}(t)=&\bar{q}_{\theta}(0)+\int_0^te^{-\mu(t-\tau)}\Big[-\frac{\theta^2}{r}s(\tau)\Big]d\tau.
\end{align*}
Thus,
\begin{align*}
  \bar{q}(t)&=\int_{\Theta}\bar{q}_{\theta}(0)dF(\theta)+\int_{\Theta}dF(\theta)\int_0^te^{-\mu(t-\tau_1)}\cr
  \cdot&\Big\{-\frac{\theta^2}{r}\int_{\tau_1}^{\infty}e^{(\rho+\mu)(\tau_1-\tau_2)}({\mathcal{A}}\bar{q})(\tau_2)
d\tau_2\Big\} d\tau_1\\
&\stackrel{\Delta}{=}(\mathcal{T}\bar{q})(t),
\end{align*}
where $({\mathcal{A}}\bar{q})(\tau_2)=\frac{1}{2}\bar{p}(0)e^{-\alpha \tau_2}+\frac{1}{2}\int_0^{\tau_2}e^{-\alpha(\tau_2-\tau_3)}(\alpha\beta-\alpha \bar{q}(\tau_3)) d\tau_3-\frac{c}{2}.$
It can be verified that $\mathcal{T}$ is a map from the Banach space $C_b([0,\infty), \mathbb{R})$ to itself. For any $\bar{q}_1,\bar{q}_2\in C_b([0,\infty), \mathbb{R})$,
\begin{align*}
&|(\mathcal{T}\bar{q}_1-\mathcal{T}\bar{q}_2)(t)|\\
&\leq \|\bar{q}_1-\bar{q}_2\|_{\infty}\int_{\Theta}\int_0^te^{-\mu(t-\tau_1)}\Big\{\frac{\theta^2}{r}\int_{\tau_1}^{\infty}e^{(\rho+\mu)(\tau_1-\tau_2)}\\
&\cdot\Big[\frac{1}{2}\int_0^{\tau_2}e^{-\alpha(\tau_2-\tau_3)}\alpha  d\tau_3\Big]
d\tau_2\Big\} d\tau_1dF(\theta)\\
&\leq  \|\bar{q}_1-\bar{q}_2\|_{\infty} \int_{\Theta}\frac{\theta^2}{2r\mu(\rho+\mu)}dF(\theta).
\end{align*}
It follows that $\mathcal{T}$ is a contraction and hence has a unique fixed point $\bar{q}\in C_b([0,\infty), \mathbb{R})$.  \hfill$\Box$

%---------------------------------------------------
\subsubsection{The case of uniform agents}
We now consider the case of uniform agents, i.e., $b_i\equiv b,$ $i=1,\cdots,N$. In this case, (\ref{eq9a})-(\ref{eq9d}) reduce to
the following equation:
\begin{eqnarray*}
\frac{d}{dt}\left[\begin{array}{c}
  \bar{p}\\\bar{q}
\\s
\end{array}\right]&=&
\left[\begin{array}{ccc}
  -\alpha &-\alpha&0\\
0&-\mu&-\frac{b^2}{r}\\
\frac{1}{2}&0 &\rho+\mu
\end{array} \right] \left[\begin{array}{c}
  \bar{p}\\\bar{q}
\\s
\end{array}\right]
+\left[\begin{array}{c}
\alpha\beta\\0\\
-\frac{c}{2}
\end{array}\right].
\end{eqnarray*}
Let
$$M=\left[\begin{array}{ccc}
  -\alpha &-\alpha&0\\
0&-\mu&-\frac{b^2}{r}\\
\frac{1}{2}&0 &\rho+\mu
\end{array} \right], \quad \bar{b}=\left[\begin{array}{c}
\alpha\beta\\0\\
-\frac{c}{2}
\end{array}\right].$$
Then
\begin{equation}\label{eq12a}
  \frac{d}{dt}\left[\begin{array}{c}
  \bar{p}\\\bar{q}
\\s
\end{array}\right]=
M\left[\begin{array}{c}
  \bar{p}\\\bar{q}
\\s
\end{array}\right]+\bar{b}.
\end{equation}
By direct computations, we have
 $$|M|={\alpha}\mu(\rho+\mu)+\frac{ \alpha b^2}{2r}$$ and the equation $Mz+\bar{b}=0$ has the solution
%$$\begin{aligned}
%  y=&\big[\frac{c}{\rho}(\alpha+\rho), -\frac{2\beta p_3}{\rho}(\alpha+\rho)^2-\frac{c}{\rho}(\alpha+\rho),\\
%&-\frac{2p_3}{\rho}(\alpha\beta p_3-\frac{c}{2}p_3+\frac{c}{2}),-\frac{c}{\rho}p_3(\alpha+\rho) \big]^T.
%\end{aligned}$$
$$\begin{aligned}
  z=&\Big[\frac{2r\beta\mu(\rho+\mu)+b^2c}{2r\mu(\rho+\mu)+b^2}, \frac{b^2(\beta-c)}{2r\mu(\rho+\mu)+b^2},\\
 &\frac{r\mu( c-\beta)}{2r\mu(\rho+\mu)+b^2}\Big]^T\hspace*{-0.1cm}.
\end{aligned}$$

Note that
 \begin{align}\label{eq17aa}
 |\lambda I-M|
=&(\lambda+\alpha)(\lambda+\mu)[\lambda-(\rho+\mu)]-\frac{\alpha b^2}{2r}\cr
=& \lambda^3+(\alpha-\rho)\lambda^2-(\mu^2+\rho\mu+\alpha\rho)\lambda\cr
&-\alpha\mu(\rho+\mu)-\frac{\alpha b^2}{2r}.
\end{align}
%Let $\lambda=\eta+\frac{\rho}{2}$. Then
%\begin{eqnarray}
%  &&|\lambda I-M|\cr
%  &=&\big(\eta-(\frac{\rho}{2}+\alpha)\big)\Big[(\eta+\frac{\rho}{2})(\eta+\frac{\rho}{2}+\alpha)(\eta-\frac{\rho}{2})-\frac{\alpha}{2}\Big]\cr
% % &=&\big(\eta-(\frac{\rho}{2}+\alpha)\big) \Big[\eta^3+(\alpha+\frac{\rho}{2})\eta^2-\frac{\rho^2}{4}\eta-(\frac{\rho^3}{8}+\frac{\alpha\rho^2}{4}+\frac{\alpha}{2})\Big]\cr
%  &\stackrel{\Delta}{=}&\big(\eta-(\frac{\rho}{2}+\alpha)\big)h(\eta).
%\end{eqnarray}
%The first column of Routh array for $h(\eta)$ is $$-1,\ -(\alpha+\frac{\rho}{2}),\  -\frac{\alpha}{2\alpha+\rho},\ \frac{\rho^3}{8}+\frac{\alpha\rho^2}{4}+\frac{\alpha}{2}.$$ By Routh's stability criterion, $h(\eta)=0$ has a root with positive real part, that is
%$h(\eta)=0$ has two roots with negative real parts. Since  $\lambda=\eta+\frac{\rho}{2}$, we have that
% $|\lambda I-M|=0$ have two roots with real parts being greater than $\rho/2$ and two roots with real parts being less than $\rho/2$.
In what follows, we use Routh's stability criterion \cite{dorf1998modern} to determine the number of roots of $|\lambda I-M|=0$ with negative real parts. The first column of the Routh array for $|\lambda I-M|$ is $-1,\ \rho-\alpha,\ \mu^2+\rho\mu+\alpha\rho+\frac{2r\alpha\mu(\rho+\mu)+\alpha b^2}{2r(\rho-\alpha)},\ \alpha\mu(\rho+\mu)+\frac{\alpha b^2}{2r}.$ It can be verified that the first column of the Routh array always has a sign change.
By Routh's stability criterion, %\cite{dorf1998modern},
(\ref{eq17aa}) has a root with a positive real part, and two roots with negative real parts.

Let $\lambda_1,\lambda_2$ be two roots of (\ref{eq17aa}) with negative real parts, and $\xi_1, \xi_2$ be the corresponding (generalized) complex eigenvectors.
Let $[\tilde{\bar{p}},\tilde{\bar{q}},\tilde{s}]^T=[\bar{p},\bar{q},s]^T-z$. The solution to equation (\ref{eq12a}) given by $z+e^{Mt}[\tilde{\bar{p}}(0),\tilde{\bar{q}}(0),\tilde{s}(0)]^T$ is in $C_{b}([0,\infty), \mathbb{R}^3)$ if and only if there exist constants $a_1, a_2$ such that $[\tilde{\bar{p}}(0),\tilde{\bar{q}}(0),\tilde{s}(0)]^T=a_1\xi_1+a_2\xi_2$.
Indeed, suppose
$$[\tilde{\bar{p}}(0),\tilde{\bar{q}}(0),\tilde{s}(0)]^T=a_1\xi_1+a_2\xi_2+a_3\xi_3,$$ where $\lambda_3$ is a root of (\ref{eq17aa}) with a positive real part and $\xi_3$ is the corresponding complex eigenvector. The solution
$$z+e^{Mt}[\tilde{\bar{p}}(0),\tilde{\bar{q}}(0),\tilde{s}(0)]^T=z+ \sum_{i=1}^3h_i(t)e^{\lambda_it}\xi_i$$
 is in $C_{b}([0,\infty), \mathbb{R}^3)$ if and only if $a_3=0$, where $h_i(t)$ are polynomials of $t$.

Denote $\xi_1=[(\xi_1^{\dag})^T,\xi_1^{\ddag}]^T$, $\xi_2=[(\xi_2^{\dag})^T,\xi_2^{\ddag}]^T$ and $z=[(z^{\dag})^T,z^{\ddag}]^T$, where
$\xi_1^{\ddag}, \xi_2^{\ddag}, z^{\ddag}\in \mathbb{R}$. Then we have
\begin{equation}\label{eq14}
  a_1\xi_1^{\dag}+a_2\xi_2^{\dag}=[\tilde{\bar{p}}(0),\tilde{\bar{q}}(0)]^T.
  \end{equation} Note that $[\tilde{\bar{p}}(0),\tilde{\bar{q}}(0)]^T=[p_0,q_0 ]^T-z^{\dag}$ is given. There exists a unique solution $(a_1, a_2)$ to (\ref{eq14}) if and only if
$\xi_1^{\dag}$ and $\xi_2^{\dag}$ are linearly independent.

From the analysis above, we have the following result.
\begin{proposition}
 (\ref{eq12a}) admits a unique solution $(s,\bar{q}_{\theta})$ such that $s$ and $\bar{q}_{\theta}$ are in $C_{b}([0,\infty), \mathbb{R})$ if and only if $\xi_1^{\dag}$ and $\xi_2^{\dag}$ are linearly independent. In this case, \textbf{A4)} holds. $\hfill \Box$
\end{proposition}

\begin{example}
 Take parameters as $[\alpha\ \beta\ \mu\ b_i\ \sigma \ \rho\ r]=[1\ 10\ 0.15\ 1\ 0.2\ 0.6\ 1]$.
Let $p(0)=1$,\ $q_i(0)\sim N(2,0.2)$. In this case,  $M$ has only two eigenvalues with negative real parts $-0.6875+0.3944i$ and $-0.6875-0.3944i$. The corresponding eigenvectors are
 $[ -0.8557, \ 0.2674 + 0.3375i,\  0.2768 + 0.0759i]^T$ and $[-0.8557 ,\  0.2674 - 0.3375i, \   0.2768 - 0.0759i]^T$, respectively. By (\ref{eq14}), we have  $a_1=  1.4429 + 7.8552i$ and $a_2= 1.4429 - 7.8552i$.  Then (\ref{eq12a}) admits a unique solution in $C_{b}([0,\infty), \mathbb{R}^3)$. However,
$\frac{b^2}{2r\mu(\rho+\mu)}=4.444>1$. The parameters in this example satisfy the condition of Proposition 2, but not of Proposition 1.
\end{example}

\subsection{$\varepsilon$-Nash equilibrium}

Consider the system of $N$ firms. Let the control strategy of firm $i$ be given by
\begin{equation}\label{eq16a}
\hat{u}_i=-\frac{b_i}{r}s, \quad  i=1,\cdots,N,
\end{equation}
where $s\in C_{b}([0,\infty), \mathbb{R})$ is determined by the equation system (\ref{eq9a})-(\ref{eq9d}).
After the strategy (\ref{eq16a}) is applied, the closed-loop dynamics for firm $i$ may be written as follows:
\begin{eqnarray}\label{eq17a}
\frac{d\hat{p}(t)}{dt}&=&-\alpha \hat{p}(t)-\alpha \hat{q}^{(N)}(t)+\alpha\beta,\\ \label{eq18a}
d\hat{q}_i(t)&=&-\mu \hat{q}_i(t)dt-\frac{b_i^2}{r}s(t)dt+\sigma dw_i,\ i=1,\cdots,N.
\end{eqnarray}
Denote $\varepsilon_N=
     \Big|\int_{\Theta}\theta^2dF_N(\theta)- \int_{\Theta}\theta^2dF(\theta)\Big|$.
\begin{theorem}\label{thm2}
  For the system (\ref{eq1})-(\ref{eq2}), if assumptions \textbf{A1)}-\textbf{A4)} hold, then the closed-loop system (\ref{eq17a})-(\ref{eq18a}) satisfies
      \begin{align}\label{eq19}
     &\sup_{t\geq 0, N\geq 1 }\big\{E|\hat{p}(t)|^2+E|\hat{q}^{(N)}(t)|^2\big\}\leq C_0,\\
     \label{eq20a}
    &\sup_{t\geq 0}E \big\{ |\hat{p}(t)-\bar{p}(t)|^2+|\hat{q}^{(N)}(t)-\bar{q}(t)|^2\big\}\cr
 =&O(\varepsilon_N^2+{1}/{N}).
     \end{align}
    % where $\bar{p}^*$ is the solution to the first equation in (\ref{eq4}) when $f$ is replaced by $f^*$.
\end{theorem}

\emph{Proof.} By (\ref{eq18a}), it follows that
\begin{eqnarray*}\label{eq22a}
  d\hat{q}^{(N)}(t)&=&\big[-\mu\hat{q}^{(N)}(t)- \frac{1}{N}\sum_{i=1}^N\frac{b_i^2}{r}s(t)\big]dt
+\frac{1}{N}\sum_{i=1}^N\sigma dw_i(t).
\end{eqnarray*}
From this together with (\ref{eq17a}), we have
\begin{align}\label{eq22aa}
& \left[\begin{array}{c}
 d\hat{p}\\d\hat{q}^{(N)}
\end{array}\right]=\left[\begin{array}{cc}
-\alpha&-\alpha \\
0&-\mu
\end{array}\right]\left[\begin{array}{c}
 \hat{p}\\\hat{q}^{(N)}
\end{array}\right]dt\cr
&-
\left[\begin{array}{c}
\alpha\beta\\ \frac{1}{N}\sum_{i=1}^N\frac{b_i^2}{r}s
\end{array}\right]dt+\left[\begin{array}{c}
0\\\frac{1}{N}\sum_{i=1}^N\sigma dw_i(t)
\end{array}\right].
\end{align}
%It can be verified that $\left[\begin{array}{cc}
%-\alpha&-\alpha \\
%\frac{1}{N}\sum_{i=1}^N \frac{b_i^2}{2r(\rho+\alpha)}&-\mu
%\end{array}\right]$ is Hurwitz.
By $s\in C_{b}([0,\infty),\mathbb{R})$ and elementary linear SDE estimates, we have
$$\sup_{t\geq 0, N\geq 1}\big\{E|\hat{p}(t)|^2+E|\hat{q}^{(N)}(t)|^2\big\}\leq C_0.$$

Denote $\eta=[\hat{p}-\bar{p}, \hat{q}^{(N)}-\bar{q}]^T$. By (\ref{eq9a}), (\ref{eq9c}) and (\ref{eq22aa}), we have
$$
 d\eta={G}\eta dt+\left[\begin{array}{c}
 0\\
\Delta_s
\end{array}\right]dt+\left[\begin{array}{c}
0\\\frac{1}{N}\sum_{i=1}^N\sigma dw_i
\end{array}\right],$$
where
\begin{equation}\label{eq22bb}
{G}=\left[\begin{array}{cc}
-\alpha&-\alpha \\
0&-\mu
\end{array}\right],
\end{equation}
$$\Delta_s=\int_{\Theta}\frac{\theta^2}{r}sdF(\theta)-\frac{1}{N}\sum_{i=1}^N\frac{b_i^2}{r}s.$$
%Thus,
%\begin{equation}\label{eq24a}
%  \eta(t)=e^{{G}t}\eta(0)+\int_0^te^{{G}(t-\tau)}\left[\begin{array}{c}
% 0\\
%\Delta_s(\tau)
%\end{array}\right]d\tau+\zeta(t),
%  \end{equation}
%where $$\zeta(t)=\int_0^te^{{G}(t-\tau)}\left[\begin{array}{c}
%0\\\frac{1}{N}\sum_{i=1}^N\sigma dw_i(\tau)
%\end{array}\right].$$
By solving this linear SDE and using the fact that ${G}$ is Hurwitz,
%As $N\to \infty$, we have
%$$\begin{aligned}
%\Delta_s(t)=&O(\varepsilon_N),\\
%E\|\zeta(t)\|^2\leq&\frac{\sigma^2}{N} \int_0^t\|e^{G(t-\tau)}\|^2d\tau=O(\frac{1}{N}),
%\end{aligned}$$
%and
%$$\begin{aligned}
%  E\|\eta(0)\|^2=&E\|q^{(N)}(0)-q_0\|^2\\
%  \leq&\frac{\max_{1\leq i\leq N}Eq_i^2(0)-q_0^2}{N}=O(\frac{1}{N}).
%\end{aligned}$$
%Hence, it follows from (\ref{eq24a}) that
we can show $$ \sup_{t\geq 0}E\|\eta(t)\|^2=O(\varepsilon^2_N+\frac{1}{N}),$$ which leads to (\ref{eq20a}).  \hfill $\Box$

By the above theorem, we can obtain the next corollary.
\begin{corollary}
  For the system (\ref{eq1})-(\ref{eq2}), if assumptions \textbf{A1)}-\textbf{A4)} hold, then the closed-loop system (\ref{eq17a})-(\ref{eq18a}) satisfies
  \begin{align*}
   &   \sup_{t\geq 0, N\geq 1 }E\int_0^{\infty}e^{-\rho t}\big\{|\hat{p}(t)|^2+|\hat{q}^{(N)}(t)|^2\big\}dt\leq C_0,\\
    &E\int_0^{\infty}e^{-\rho t} \big\{ |\hat{p}(t)-\bar{p}(t)|^2+|\hat{q}^{(N)}(t)-\bar{q}(t)|^2\big\}dt\cr
  =&O(\varepsilon_N^2+{1}/{N}).
     \end{align*}
\end{corollary}

We are now in a position to show an asymptotic Nash equilibrium property. Denote
  \begin{align*}
    {\mathcal U}_c=&\Big\{u_i: u_i(t)\ \hbox{is adapted to} \ \sigma\{\cup_{j=1}^N  {\mathcal F}_t^{j}\},\cr
&E\int_0^{\infty} e^{-\rho t} u_i^2(t)dt<\infty\Big\}.\cr
\hat{u}_{-i}=&(\hat{u}_{1}, \cdots, \hat{u}_{i-1}, \hat{u}_{i+1}, \cdots, \hat{u}_{N}).
  \end{align*}
\begin{theorem}\label{thm3}
 For the problem (\ref{eq1})-(\ref{eq2}), assume that \textbf{A1)}-\textbf{A4)} hold. Then the set of strategies $(\hat{u}_1,\cdots,\hat{u}_N)$ given by (\ref{eq16a}) is an $\varepsilon$-Nash equilibrium, i.e.,
 \begin{equation}\label{eq26n}
   J_i(\hat{u}_i, \hat{u}_{-i})-\varepsilon\leq \inf_{u_i\in {\mathcal U}_c}J_i({u}_i, \hat{u}_{-i})\leq J_i(\hat{u}_i, \hat{u}_{-i}),
 \end{equation}
 where $\varepsilon=O(\varepsilon_N+\frac{1}{\sqrt{N}})$.
\end{theorem}

\emph{Proof.} See Appendix A.  \hfill $\Box$

\section{Social Solutions to Output Adjustment}

We first construct an auxiliary optimal control problem by examining the social cost variation due to the control perturbation of a single
agent. %and then show the relation between the original problem and
 %the auxiliary optimal control problem.
{Then, by mean field approximations we design a set of decentralized strategies which is shown to have asymptotic social optimality.}

\subsection{An auxiliary optimal control problem}

We now provide a property of the social optimum problem which implies that $J^{(N)}_{\rm soc}$ has a minimizer.
\begin{lemma}
  $J^{(N)}_{\rm soc}(u)$ is coercive with respect to $(u_1,\cdots, u_N)$, i.e., there exist constants $C_2>0$ and $C_3>0 $ such that
  $$J^{(N)}_{\rm soc}(u)\geq \frac{C_2}{N}E\int_0^{\infty} e^{-\rho t}\sum_{i=1}^Nu_i^2 dt-C_3.$$
\end{lemma}
\emph{Proof.} From Lemma \ref{lem3}, we can get the lemma immediately.  \rightline{$\Box$}

This lemma ensures the existence of a centralized optimal solution to the social optimum problem in (\ref{eq1})-(\ref{eq2}) and (\ref{eq4}) (see \cite{KZ2006}).

We now derive an auxiliary optimal control problem from the original social optimum problem by perturbing the strategy of a fixed agent. %and further construct an auxiliary optimal control problem to approximate
%the equivalent problem. %by person-by-person optimality. Specifically
%We provide an equivalent transformation of the original social optimum problem by analyzing open-loop centralized strategy.
Denote the control problem \textbf{(P1)}:
 \begin{align*}
  \frac{d{p}}{dt}=&-\alpha p-\frac{\alpha}{N}q_i-\alpha \hat{q}_{-i}^{(N)}+\alpha\beta,\cr
    dq_i=&-\mu q_idt+b_iu_idt+\sigma dw_i,\cr
    \frac{dv_i}{dt}=&-\alpha v_i-\alpha q_i, \ v_i(0)=0,\cr
J^{*}_i(u_i)=&E \int_0^{\infty}e^{-\rho t}\Big[(c-p)q_i-v_i \hat{q}^{(N)}_{-i}+ru_i^2\Big]dt,
     \end{align*}
     where $\hat{q}_{-i}^{(N)}=\frac{1}{N}\sum_{j=1,j\neq i}^N \hat{q}_j$ and $\hat{u}_{-i}$ is given.

\begin{lemma}
  If $\hat{u}=(\hat u_1, \ldots, \hat u_N)$ minimizes $J^{(N)}_{\rm soc}$ where each $\hat{u}_i\in {\mathcal U}_{c}$, then $\hat{u}_i$ is
  necessarily the optimal strategy of Problem \textbf{(P1)}.
\end{lemma}
{\emph{Proof.}} It follows from (\ref{eq3}) that
\begin{align*}
J^{(N)}_{\rm soc}=&\frac{1}{N} \sum_{i=1}^N J_i(\hat{u}_1,\cdots,\hat{u}_{i-1},u_i,\hat{u}_{i+1},\cdots,\hat{u}_{N})\cr
 %= E \int_0^{\infty}e^{-\rho t}[(c-{p}) \sum_{j=1}^N{q}_j+\sum_{j=1}^Nr {u}_j^2 ]
  =&\frac{1}{N}E \int_0^{\infty} e^{-\rho t}(Y_i+Y_i^{\prime})dt,
\end{align*}
where $$Y_i=( c-p){q}_i-{p} \sum_{j=1, j\not= i}^N\hat{q}_j+r {u}_i^2,$$
$$Y_i^{\prime}= c\sum_{j=1, j\not= i}^N\hat{q}_j+\sum_{j=1, j\not = i}^Nr \hat{u}_j^2.$$
By (\ref{eq1}),
\begin{align*}
p(t)%&=e^{-\alpha t}p(0)+\int_0^{t}e^{-\alpha(t-\tau)}(\alpha\beta-\alpha q^{(N)})d\tau\cr
=&e^{-\alpha t}p(0)+\int_0^{t}e^{-\alpha(t-\tau)}(\alpha\beta-\alpha \hat{q}^{(N)}_{-i})d\tau\cr
&-\frac{\alpha}{N}\int_0^{t}e^{-\alpha(t-\tau)}q_id\tau.
\end{align*}
 Thus, $Y_i=Z_i+Z_i^{\prime}$
where
$$Z_i=(c-p)q_i+{\alpha}\int_0^{t}e^{-\alpha(t-\tau)}q_id\tau\cdot \hat{q}^{(N)}_{-i}+ru_i^2,$$
$$Z_i^{\prime}=-\Big[e^{-\alpha t}p(0)+\int_0^{t}e^{-\alpha(t-\tau)}(\alpha\beta-\alpha \hat{q}^{(N)}_{-i})d\tau\Big]\sum_{j=1, j\not= i}^N\hat{q}_j.$$
Note that
$J^{(N)}_{\rm soc}=\frac{1}{N}E \int_0^{\infty} e^{-\rho t}(Z_i+Z_i^{\prime}+Y_i^{\prime})dt,$ where neither $Z_i^{\prime}$ nor $Y_i^{\prime}$ changes with $u_i$.
Thus, for agent $i$ minimizing $J^{(N)}_{\rm soc}$ is equivalent to minimizing
 $E\int_0^{\infty} e^{-\rho t} Z_i dt$,  which in turn is equal to $J_i^*(u_i)$.
%\begin{align*}
%  &J^{*}_i(u_i)\cr
%  =&E \int_0^{\infty}e^{-\rho t}\Big[(c-p)q_i+{\alpha}\int_0^{t}e^{-\alpha(t-\tau)}q_id\tau\cdot \hat{q}^{(N)}_{-i}+ru_i^2\Big]dt.
%  \end{align*}
%Let $v_i=-\alpha\int_0^{t}e^{-\alpha(t-\tau)}q_id\tau$. Then
%\begin{align*}
%  J^{*}_i(u_i)&=E \int_0^{\infty}e^{-\rho t}\Big[(c-p)q_i-v_i \hat{q}^{(N)}_{-i}+ru_i^2\Big]dt,\\
%\frac{dv_i}{dt}&=-\alpha v_i-\alpha q_i, \ v_i(0)=0.
%\end{align*}
\hfill$\Box$

\subsection{Mean field approximation}

To approximate Problem (P1) for large $N$, % to approximate J^{*}_i %$\sum_{j=1}^{N}J_j$,
 we construct the auxiliary limiting optimal control problem \textbf{(P2)}:
\begin{eqnarray}\label{eq5s}
&&\left[\begin{array}{c}
  d\bar{p}\\dq_i\\dv_i
\end{array}\right]=
\left[\begin{array}{ccc}
 -\alpha&0&0\\0&-\mu&0\\
 0&-\alpha&-\alpha
\end{array}\right]
\left[\begin{array}{c}
  \bar{p}\\q_i\\v_i
\end{array}\right]dt+
\left[\begin{array}{c}
  0\\b_i\\0
\end{array}\right]u_idt\cr
&&+\left[\begin{array}{c}
 \alpha\beta-\alpha \bar{q} \\0\\0
\end{array}\right]dt+\left[\begin{array}{c}
  0\\1\\0
\end{array}\right]dw_i, \ \left[\begin{array}{c}
  \bar{p}(0)\\q_i(0)\\v_i(0)
\end{array}\right]=  \left[\begin{array}{c}
p_0\\q_i(0)\\0
\end{array}\right]
\end{eqnarray}
with cost function
\begin{align}\label{eq6s}
\bar{J}_i^{*}(u_i)=E\int_0^{\infty} e^{-\rho t}[(c-{\bar{p}}) {q}_i-\bar{q} v_i+r {u}_i^2]dt.
\end{align}
Here $\bar{q}\in C_{\rho/2}([0,\infty), \mathbb{R})$ is a deterministic function, which is an approximation of $q^{(N)}$ for large $N$.

For the system (\ref{eq5s})-(\ref{eq6s}),  we take $\bar{p}$ as an exogenous signal. %It is easy to get that $ \bar{J}^*_i(u_i)$ is strictly convex.

\begin{lemma}\label{lem1b}
  $ \bar{J}^*_i(u_i)$ is strictly convex and coercive.
\end{lemma}
{\emph{Proof.}}  {For the system (\ref{eq5s})-(\ref{eq6s}), the state is $(q_i,v_i)$; $\bar{p}$ and $\bar{q}$ are
 not dependent on the control $u_i$.
%By a direct method for verifying the definition of convexity, w
We can directly show that
$\bar{J}^*_i(u_i)$ is strictly convex in $u_i$.
%It is straightforward to verify $J^{*\prime}_i(u_i)$ is positive definite.
 Following the proof of Proposition \ref{prop1}, we can show that  $ \bar{J}^*_i(u_i)$ is coercive.}
%  It is a direct observation that there exists a small constant $\epsilon_0>0$ such that
%$$\bar{J}_{i,\epsilon_0}(u_i)\stackrel{\Delta}{=}\bar{J}_i(u_i)-2\epsilon_0\|u_i\|^2_{\rho} $$ is convex. Notice that
% $\left[\begin{array}{cc}
% -\alpha-\frac{\rho}{2}&0\\0&-\frac{\rho}{2}
%\end{array}\right]$ is Hurwitz. By Lemma A. 1 in \cite{huang2010large}, there exists a constant $C$ dependent of $p_0$ and $q_0$ such that $\sup_{\|u_i\|_{\rho}\leq 1}\big|\bar{J}_{i,\epsilon_0}(u_i)\big|\leq C$. %It can be verified that $\bar{J}_{i,\epsilon_0}(u_i)$ is convex since $\epsilon_0$ is small.
%For $\|u_i\|_{\rho}\geq 1$, it follows from the convexity of  $\bar{J}_{i,\epsilon_0}(u_i)$ that
%$$\begin{aligned}
%\bar{J}_{i,\epsilon_0}\Big(\frac{u_i}{\|u_i\|_{\rho}}\Big)&\leq \frac{1}{\|u_i\|_{\rho}}\bar{J}_{i,\epsilon_0}({u_i})+\frac{\|u_i\|_{\rho}-1}{\|u_i\|_{\rho}}\bar{J}_{i,\epsilon_0}({0})\\
%&\leq
%\frac{1}{\|u_i\|_{\rho}}\bar{J}_{i,\epsilon_0}({u_i})+C,
%\end{aligned}$$
%which implies $\bar{J}_{i,\epsilon_0}({u_i})\geq -2C \|u_i\|_{\rho}.$
%Hence for any $u_i$, $$\bar{J}_{i,\epsilon_0}({u_i})\geq -C(2\|u_i\|_{\rho}+1),$$
%which further leads to
%$$\bar{J}_i(u_i)\geq \epsilon_0 \|u_i\|_{\rho}^2 -C.$$
%\rightline{$\Box$}
 \hfill $\Box$

    Let $$A=\left[
\begin{array}{cc}
-\mu&0\\
-\alpha&-\alpha
\end{array}
\right],\quad B_i=\left[
\begin{array}{c}b_i\\
0\end{array}
\right].$$
By Theorem \ref{thm1a}, Lemma \ref{lem1b} and \cite{KZ2006}, Problem (P2) has the unique optimal control given by $$u_i=-\frac{1}{r}B_i^T\check{s},$$
where $\check{s}\in {C}_{\rho/2}([0,\infty), \mathbb{R}^2)$ is determined from
  $$ \rho \check{s} =\frac{d\check{s}}{dt}+A^T \check{s}+
  \frac{1}{2}[c-\bar{p}, -\bar{q}]^T.$$
  Let $B_{\theta}=[\theta,0]^T$. When $\theta=b_i$, we have $B_{\theta}=B_{i}$, and $y_{\theta}=y_i$. $y_{\theta}$ is regarded as the expectation of $(q_{\theta},v_{\theta})^T$ given the parameter ${\theta}$ in the individual dynamics.
Following the standard approach in mean field control \cite{HCM07}, \cite{HCM12}, we construct the equation system as follows:
\begin{align}\label{eq9}
\frac{d\bar{p}}{dt}=&\alpha [-\bar{p}+\beta-\bar{q}],\\ \label{eq9e}
  \rho \check{s} =&\frac{d\check{s}}{dt}+A^T \check{s}+
  \frac{1}{2}[c-\bar{p}, -\bar{q}]^T,\\ \label{eq9f}
\frac{dy_{\theta}}{dt}=&A y_{\theta}-\frac{1}{r}B_{\theta}B^T_{\theta}\check{s},\\ \label{eq9g}
\bar{q}=&[1,0]\int_{\Theta}y_{\theta}dF(\theta),
\end{align}
where $y_{\theta}(0)
=[q_0,0]^T$.
For further analysis, we assume

\textbf{A5)} There exists a solution $(\check{s}, y_{\theta}, \theta \in \Theta)$ to (\ref{eq9})-(\ref{eq9g}) such that for any $\theta \in \Theta$, both $\check{s}$ and $y_{\theta}$ are within $C_b([0,\infty), \mathbb{R}^2)$. % which is bounded and continuous.

 For the case of uniform agents $(b_i\equiv b)$, the equation system (\ref{eq9})-(\ref{eq9g}) reduces to
\begin{align}\label{eq10c}
\frac{d\bar{p}}{dt}=&\alpha [\beta-\bar{p}-\bar{q}],\\ \label{eq10e}
  \rho \check{s} =&\frac{d\check{s}}{dt}+A^T \check{s}+
  \frac{1}{2}[c-\bar{p}, -\bar{q}]^T,\\ \label{eq10f}
\frac{dy}{dt}=&A y-\frac{1}{r}BB^T\check{s},\
\end{align}
where $\bar{p}(0)=p_0$, $y(0)=[q_0,0]^T$ and $B=[b, 0]^T$.
Let
\begin{equation} \label{eq10d}
  M_s=\left[
\begin{array}{ccccc}
  -\alpha&0&0&-\alpha&0\\
  \frac{1}{2}&\rho+\mu&\alpha&0&0\\
  0&0&\rho+\alpha& \frac{1}{2}&0\\
  0&-\frac{b^2}{r}&0&-\mu&0\\
  0&0&0&-\alpha&-\alpha
\end{array}\right],\
 \end{equation}
$$ \bar{b}_s=\left[
 \begin{array}{c}
 \alpha\beta\\
 -\frac{c}{2}\\
   0\\0\\0
 \end{array}\right], \qquad \varphi=\left[
 \begin{array}{c}
 \bar{p}\\
\check{s}\\
   y
 \end{array}\right].
$$
 Then
\begin{equation}\label{eq13n}
\frac{d \varphi}{dt}=M_s\varphi+\bar{b}_s.
\end{equation}
 By straightforward computation, we can show
 %$$|M_s|=-\alpha^2\Big[\mu(\rho+\alpha)(\rho+\mu)+\frac{b^2}{2r}(\rho+2\alpha)\Big],$$ and
  $M_s\varphi+\bar{b}_s=0$ has a unique solution, denoted as $z_s$.
%$$\begin{aligned}
%  z_s=&\Big[\frac{2r\beta\mu(\rho+\mu)+b^2c}{2r\mu(\rho+\mu)+b^2}, \frac{r\mu(\beta-c)}{b^2+2r\mu(\rho+\mu)},\\
%  &\frac{(\beta-c)(2\alpha r\mu-b^2 )}{2(\alpha+\rho)[2r\mu(\rho+\mu)+b^2]},\\
%  & \frac{(\beta-c)b^2}{b^2+2r\mu(\rho+\mu)},-\frac{(\beta-c)b^2}{2r\mu(\rho+\mu)+b^2} \Big]^T.
%\end{aligned}$$
Furthermore, we have
$$\begin{aligned}
   & |\lambda I-M_s|
=(\lambda+\alpha)\times\Big[ \lambda^4-2\rho\lambda^3\\
&+(\rho^2-(\alpha+\mu)\rho-\alpha^2-\mu^2)\lambda^2+\rho[(\alpha+\mu)\rho+\alpha^2+\mu^2]\lambda\\
  &+\alpha\mu(\rho+\alpha)(\rho+\mu)+\frac{\alpha b^2}{2r}(\rho+2\alpha)\Big].\\
\end{aligned}$$
%Let
%$$\begin{aligned}
%g(\lambda)=&\lambda^3+(\alpha-\rho)\lambda^2-(\mu^2+\rho\mu+\alpha\rho)\lambda\\
%  &-\alpha\mu(\rho+\mu)+\frac{\alpha b^2}{2r}.
%\end{aligned}$$
%We now use Routh's stability criterion \cite{dorf1998modern} to determine the number of roots of $g(\lambda)$ with negative real parts. The first column of the Routh array for $g(\lambda)$ is $1,\ -2\rho,b_1,\rho(a_3+\frac{2\rho a_4}{b_1}), a_4,$ where $b_1=\rho^2-\frac{\alpha+\mu}{2}\rho-\frac{\alpha^2+\mu^2}{2}$. It can be verified that the first column of the Routh array always has two sign changes.
By Routh's stability criterion \cite{dorf1998modern}, we obtain that $|\lambda I-M_s|=0$ has
two roots with positive real parts, and three roots with negative real parts.

 Let $\check{\lambda}_1,\check{\lambda}_2, \check{\lambda}_3$ be two roots of $|\lambda I-M_s|=0$ with negative real parts, and $\zeta_1, \zeta_2, \zeta_3$ be the corresponding (generalized) complex eigenvectors.
Let $\tilde{\varphi}(0)= [\tilde{\bar{p}}^T(0),\tilde{\check{s}}^T(0), \tilde{y}^T(0)]^T=\varphi(0)-z_s$. The solution to equation (\ref{eq13n}) given by $z_s+e^{M_st}\tilde{\varphi}(0)$ is in $C_{b}([0,\infty), \mathbb{R}^5)$ if and only if there exist constants $\check{a}_1, \check{a}_2, \check{a}_3$ such that $\tilde{\varphi}(0)=\check{a}_1\zeta_1+\check{a}_2\zeta_2+\check{a}_3\zeta_3$.

Denote $\zeta_i=[(\zeta_i^{\dag})^T,(\zeta_i^{\ddag})^T, (\zeta_i^{\S})^T]^T$ and $\vartheta_i=[(\zeta^{\dag}_i)^T,(\zeta^{\S}_i)^T]^T$, where $\zeta_i^{\dag}\in \mathbb{R}$ and $\zeta_i^{\ddag}, \zeta_i^{\S}\in \mathbb{R}^2$, $i=1,2,3$. Then we have
\begin{equation}\label{eq14b}
  \check{a}_1\vartheta_1+\check{a}_2\vartheta_2+\check{a}_3\vartheta_3=[\tilde{\bar{p}}(0),\tilde{y}(0)]^T.
  \end{equation} Note that $[\tilde{\bar{p}}(0), \tilde{y}(0)]^T$%=[p_0,q_0,0 ]^T-
%  \left[\begin{array}{ccccc}
%1&0&0&0
%&0\\
%  0&0&0&1
%&0\\
%0&0&0&0
%&1
%\end{array}\right].z_s$$
is given. There exists a unique solution $(\check{a}_1, \check{a}_2,\check{a}_3)$ to (\ref{eq14b}) if and only if
$\vartheta_1$,  $\vartheta_2$ and $\vartheta_3$ are linearly independent.

From the analysis above, we have the following result.
\begin{proposition}
  \textbf{A5)} holds if (\ref{eq14b}) admits a solution. In particular, (\ref{eq10c})-(\ref{eq10f}) admits a unique solution $(\check{s}, y_{\theta})$ such that $\check{s}$ and $y_{\theta}$ are within $C_b([0,\infty), \mathbb{R}^2)$ if and only if $\vartheta_1$,  $\vartheta_2$ and $\vartheta_3$ are linearly independent.  $\hfill \Box$
\end{proposition}

\subsection{Asymptotic optimality}

Consider the system of $N$ firms. Let the control strategy of firm $i$ be given by
\begin{equation}\label{eq16}
\check{u}_i=-\frac{1}{r}B_i^T\check{s}, \quad  i=1,\cdots,N,
\end{equation}
where $\check{s}\in C_b([0,\infty), \mathbb{R}^2) $ is determined by the equation system (\ref{eq9})-(\ref{eq9g}).
After the strategy (\ref{eq16}) is applied, the closed-loop dynamics for firm $i$ may be written as follows:
\begin{align}\label{eq17}
\frac{d\check{p}}{dt}=&-\alpha \check{p}-\alpha \check{q}^{(N)}+\alpha\beta,\\ \label{eq18}
d\check{q}_i=&-\mu \check{q}_i-\Big[\frac{b_i^2}{r},0\Big]\check{s}dt+\sigma dw_i,\\ \label{eq19}
\frac{d\check{v}_i}{dt}=&-\alpha \check{q}_i-\alpha \check{v}_i,
\end{align}
where $\check{p}(0)=p_0$, $\check{q}_i(0)=q_i(0)$, and $\check{v}_i(0)=0$.

\begin{theorem} \label{thm4}
  Assume that \textbf{A1)}-\textbf{A3)} and \textbf{A5)} hold. The set of strategies $\{\check{u}_i=-\frac{1}{r}B_i^T\check{s},i=1,\cdots,N\}$ has asymptotic social optimality, i.e.,
  $$\left|J^{(N)}_{\rm soc}(\check{u})-\inf_{u_i\in \mathcal{U}_c}J^{(N)}_{\rm soc}({u})\right|=O(\frac{1}{\sqrt{N}}+\epsilon_N). $$
\end{theorem}
\emph{Proof.} See Appendix B.  $\hfill \Box$

We now give a closed-form expression of the asymptotic social cost.
\begin{theorem} \label{thm5}
 Assume that \textbf{A1)}-\textbf{A3)} and \textbf{A5)} hold. Then the asymptotic optimal social cost is given by
 $$\begin{aligned}
 \lim_{N\to \infty} \inf_{u_i\in {\mathcal U}_c}J^{(N)}_{\rm soc}(u)=&2\check{s}_1(0)q_0+\int_{\mathbb{R}}\check{g}_{\theta}(0)dF(\theta)\cr
 &+\int_0^{\infty}e^{-\rho t}\bar{q}(t)\bar{v}(t) dt,
  \end{aligned}$$
  where $\check{s}_1\stackrel{\Delta}{=}[1,0]\check{s}$, $\bar{v}(t)=-\alpha\int_0^te^{-\alpha(t-\tau)}\bar{q}(\tau)d\tau$ and $\check{g}_i\in C_{b}([0,\infty), \mathbb{R})$ satisfies
  $$\rho \check{g}_i=\frac{d\check{g}_i}{dt}-\frac{b_i^2}{r}\check{s}_{1}^2.$$
 \end{theorem}
 \emph{Proof.} From Lemma \ref{lem1} and Schwarz's inequality, it follows that
 \begin{equation*}
  \max_{1\leq i\leq N}\left|J_i(\check{u})-E\int_0^{\infty}e^{-\rho t}[(c-\bar{p})\check{q}_i+r\check{u}_i^2]dt\right|=O(\varepsilon_N+\frac{1}{\sqrt{N}}).
  \end{equation*}
Let $J_i^*(\check{u}_i)$ be the optimal cost of Problem (P2).
  By Theorem \ref{thm1a}, we have
\begin{align}\label{eq35}
  &\frac{1}{N}\sum_{i=1}^NE\int_0^{\infty}e^{-\rho t}[(c-\bar{p})\check{q}_i+r\check{u}_i^2]dt\cr
  =&\frac{1}{N}\sum_{i=1}^N J_i^*(\check{u}_i)+\frac{1}{N}\sum_{i=1}^NE\int_0^{\infty}e^{-\rho t}\bar{q}\check{v}_idt\cr
  =& \frac{1}{N}\sum_{i=1}^N\Big\{  2\check{s}(0)[q_0,0]^T+\check{g}_i(0)\cr
  &-\alpha E\int_0^{\infty}e^{-\rho t}\bar{q}(t)\int_0^{t}e^{-\rho \tau}\check{q}_i(\tau)d\tau dt\Big\}\cr
    =&  2\check{s}_1(0)q_0+\frac{1}{N}\sum_{i=1}^N\check{g}_i(0)+ E\int_0^{\infty}e^{-\rho t}\bar{q}(t)\check{v}^{(N)}(t)dt,
\end{align}
where% $\check{s}=[\check{s}_1,\check{s}_2]^T$,
$$\check{v}_i(t)=- \alpha\int_0^{t}e^{-\rho \tau}{\check{q}}_i(\tau)d\tau, $$
and $\check{g}_i\in C_{b}([0,\infty), \mathbb{R})$ satisfies
$$\rho \check{g}_i=\frac{d\check{g}_i}{dt}-\frac{b_i^2}{r}\check{s}^2_{1}.$$
By Schwarz's inequality and Lemma \ref{lem1},
\begin{align}\label{eq36}
&\left|E\int_0^{\infty}e^{-\rho t}\bar{q}(t)|\check{v}^{(N)}(t)-\bar{v}(t)|dt\right|^2\cr
\leq& E\int_0^{\infty}e^{-\rho t}|\bar{q}(t)|^2dt \cdot E\int_0^{\infty}e^{-\rho t} | \check{v}^{(N)}(t)-\bar{v}(t)|^2dt\cr
\leq& C(\frac{1}{N}+\varepsilon^2_N).
\end{align}
Notice that $\check{g}_{\theta}(0)$ is continuous in $\theta$. By the weak convergence of $F_N$ to $F$, we obtain
$$\lim_{N\to \infty}\frac{1}{N}\sum_{i=1}^N\check{g}_i(0)=\int_{\mathbb{R}}\check{g}_{\theta}(0)dF(\theta),$$
which together with (\ref{eq35}) and (\ref{eq36}) completes the proof.  $\hfill \Box$

%\begin{remark}
%From Theorems \ref{thm1a} and \ref{thm5}, we find that formally the optimal social cost has an additional term $\int_0^{\infty}e^{-\rho t}\bar{q}(t)\bar{v}(t) dt$, compared with the optimal Nash cost. Note that $\bar{v}(t)=-\alpha\int_0^te^{-\alpha(t-\tau)}\bar{q}(\tau)d\tau<0$. The optimal social cost
%  is likely smaller than the optimal Nash cost.
%\end{remark}

\section{Numerical Simulation}

In this section, we provide a numerical example to illustrate the evolution of firms' outputs and compare prices, average outputs and optimal costs under two different solution frameworks: the mean field game and social optimization.

 Take %parameters as
  $[\alpha,\ \beta,\ \mu,\ b_i,\ \sigma, \ \rho\ r]=[1,\ 10,\ 0.15,\ 1,\ 0.2,\ 0.6\ 1]$.
Let $p(0)=1$, and $q_i(0)\sim N(2,0.2)$. %In this case, $M$ only has two eigenvalues with negative real parts $-0.6875+0.3944i$ and $-0.6875-0.3944i$.
 %The corresponding eigenvectors are
 %$[0.8899 ,  -0.2781 - 0.3509i,    -0.0053 -0.0156i,   -0.0336 -0.0790i ]^T$ and $[0.8899 ,  -0.2781+0.3509i,    -0.0053+0.0156i,   -0.0336 +0.0790i ]^T$, respectively.
 By (\ref{eq14}), we have $a_1=1.4429 + 7.8552i$ and $a_2=1.4429- 7.8552i$; by (\ref{eq14b}), we get $\check{a}_1=9$, $\check{a}_2=-3.5906-6.5046i$ and $\check{a}_3=-3.5906+6.5046i$. Thus, we can write the exact expressions of $s$ and $\check{s}$.

     {It can be verified that \textbf{A1)}-\textbf{A5)} hold. Figs. 1 and 2 show the curves of output levels of firms within the two frameworks, respectively. They remain positive, although %some outputs occasionally get close to zero due to
random fluctuations appear and there are greater fluctuations in the social optimum framework. After the transient phase, the output levels of the firms behave similarly.}

     \begin{center}
 \vskip 0.3cm $\epsfig{figure=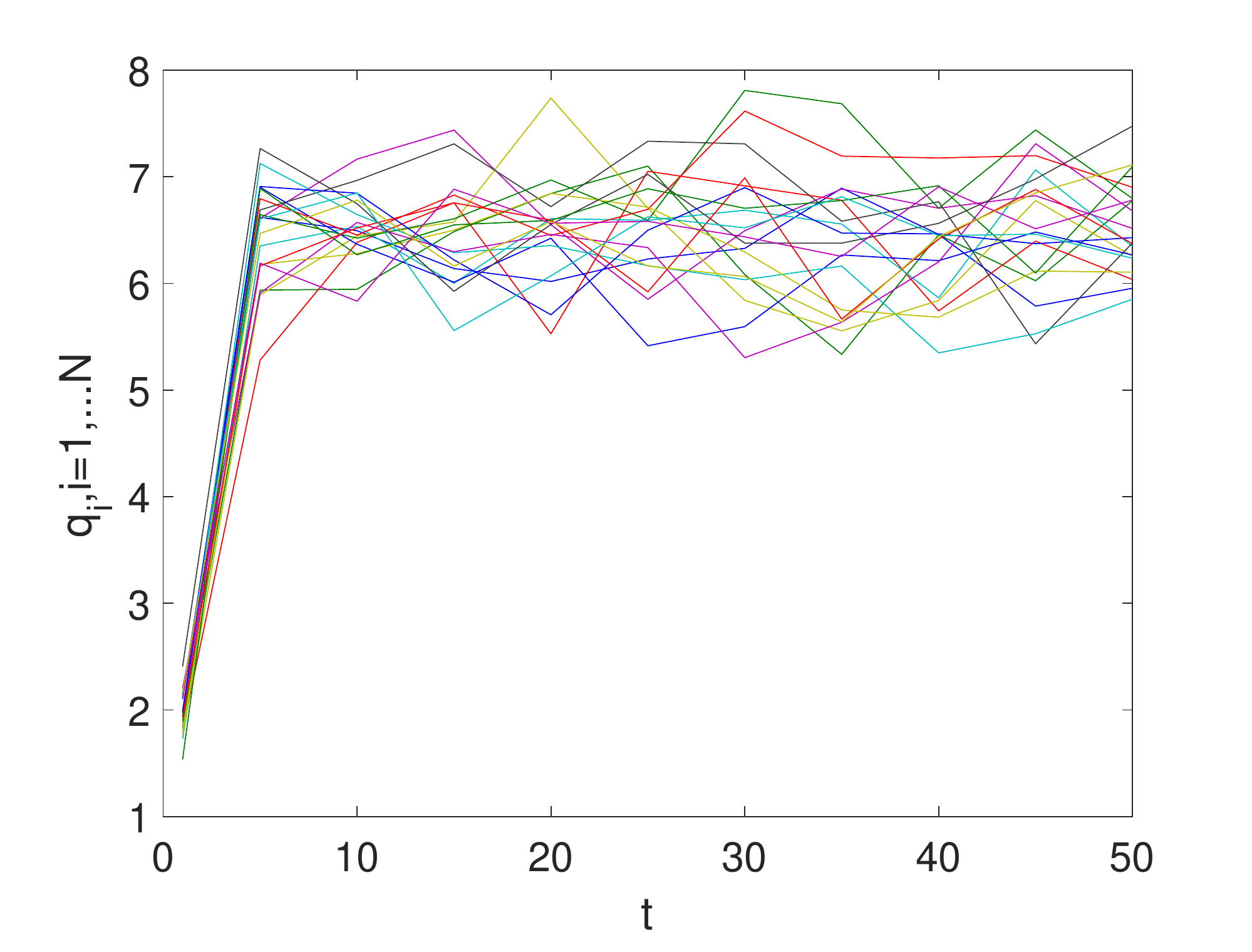,height=6cm}$\\
 {\small Fig. 1: Curves of $q_i, i=1,\cdots,20$ in the game solution}
\end{center}
 \begin{center}
 \vskip 0.3cm $\epsfig{figure=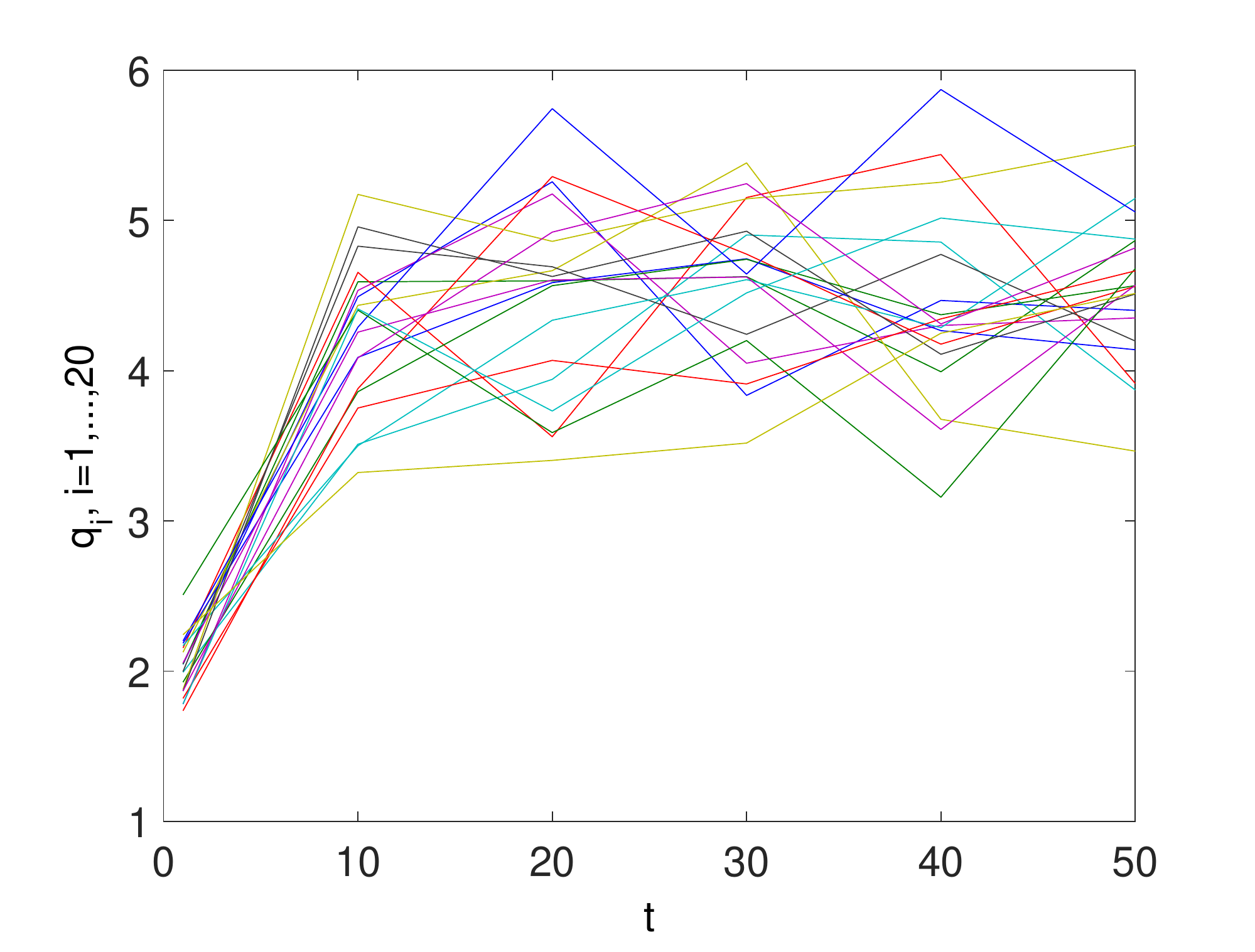,height=6cm}$\\
 {\small Fig. 2: Curves of $q_i, i=1,\cdots,20$ in the social optimum solution}
\end{center}

Fig. 3 depicts the curves of $p$ and $\bar{p}$ within game and social frameworks, when the total number of agents is 50. Fig. 4 shows the curves of $q^{(50)}$ and $q^{(\infty)}$ within two frameworks, where $q^{(\infty)}=\bar{q}$ is the average output of firms in the infinite population case.
It can be seen %from Figs. 3 and 4
 that the curves of $p$ and $\bar{p}$ as well as $q^{(50)}$ and $q^{(\infty)}$ coincide well, which illustrates the accuracy of the mean field approximation. %Compared with the game case,
 From the game framework to the social framework, the price gets a significant increase, and the average of outputs becomes lower.

     \begin{center}
 \vskip 0.3cm $\epsfig{figure=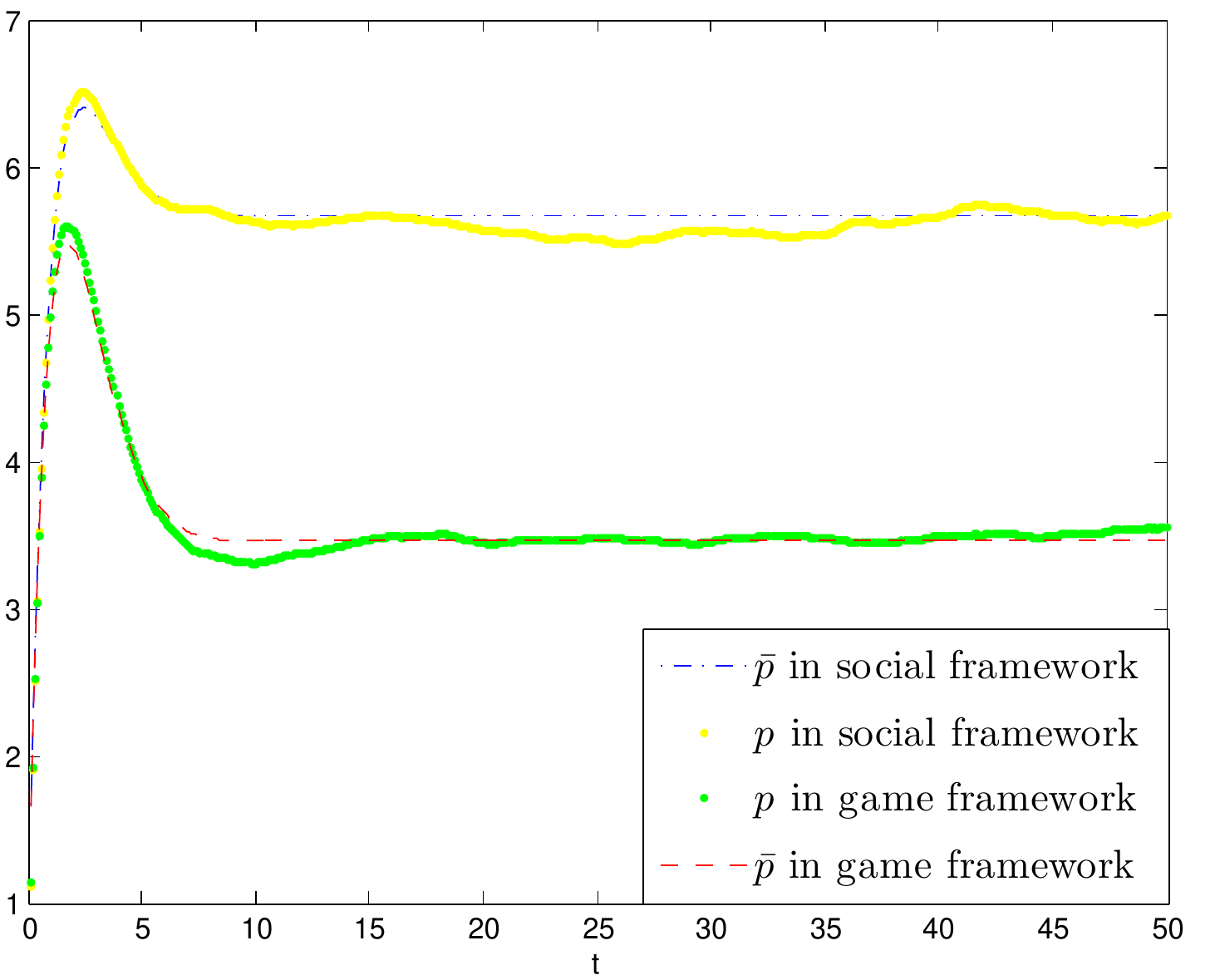,height=5.2cm}$\\
 {\small Fig. 3: Curves of $p$ and $\bar{p}$ in  the game and social optimum, respectively}
\end{center}
 \begin{center}
 \vskip 0.3cm $\epsfig{figure=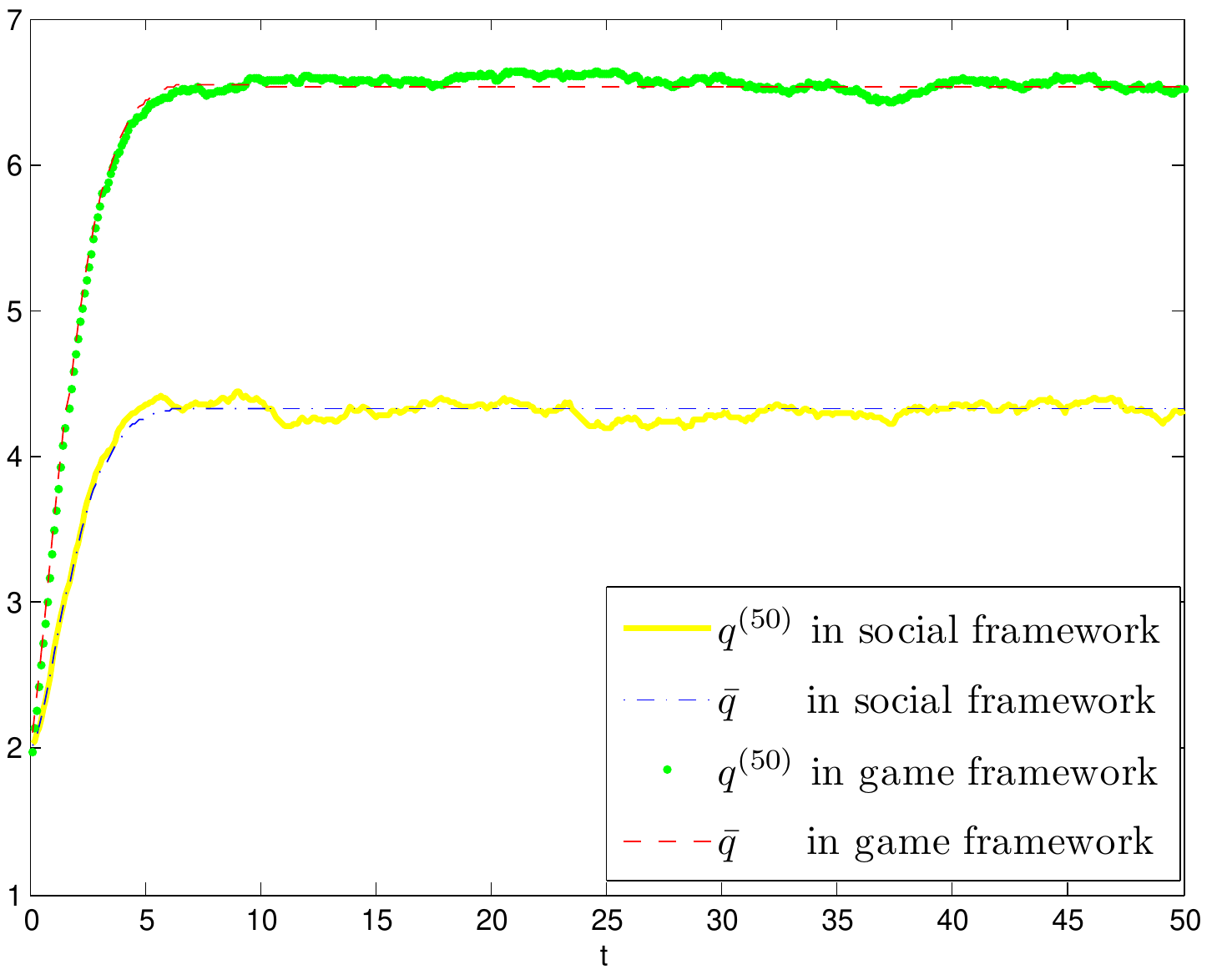,height=5.5cm}$\\
 {\small Fig. 4: Curves of $q^{(50)}$, $q^{(\infty)}$ in  the game and social optimum, respectively}
\end{center}

By Theorem \ref{thm1a}, we get that the asymptotic Nash cost is given by
$$ J_{\rm Nash}^{(\infty)}\stackrel{\Delta}{=}\lim_{N\to\infty}J_i(\hat{u}_i,\hat{u}_{-i})=2s(0)q_0+g(0),$$
where $g(0)=-\int_0^{\infty}e^{-\rho t}s^2dt.$
Note that $s=z^{\ddag}+a_1e^{\lambda_1t}\xi_1^{\ddag}+a_2e^{\lambda_2t}\xi_2^{\ddag}$.
 Then
$$\begin{aligned}
 J_{\rm Nash}^{(\infty)}=&2s(0)q_0-\frac{b^2}{r}\int_0^{\infty}e^{-\rho t}s^2(t)dt\\
=&2q_0(z^{\ddag}+a_1\xi_1^{\ddag}+a_2\xi_2^{\ddag})\cr
&-\frac{b^2}{r}\int_0^{\infty}e^{-\rho t}(z^{\ddag}+a_1e^{\lambda_1t}\xi_1^{\ddag}+a_2e^{\lambda_2t}\xi_2^{\ddag})^2dt\cr
=&2q_0(z^{\ddag}+a_1\xi_1^{\ddag}+a_2\xi_2^{\ddag})-\frac{b^2}{r}\Big[\frac{(z^{\ddag})^2}{\rho}+\frac{a_1^2(\xi_1^{\ddag})^2}{\rho-2\lambda_1}\cr
&+\frac{a_2^2(\xi_2^{\ddag})^2}{\rho-2\lambda_2}
+\frac{2a_1z^{\ddag}\xi_1^{\ddag}}{\rho-\lambda_1}+\frac{2a_2z^{\ddag}\xi_2^{\ddag}}{\rho-\lambda_2}+\frac{2a_1a_2\xi_1^{\ddag}\xi_2^{\ddag}}
{\rho-\lambda_1-\lambda_2}\Big].
\end{aligned}$$

By Theorem \ref{thm5},  $$\begin{aligned}
 J_{\rm soc}^{(\infty)}\stackrel{\Delta}{=}\lim_{N\to \infty}J^{(N)}_{\rm soc}(\check{u})=&2\check{s}_1(0)q_0+\check{g}(0)\cr
 &+\int_0^{\infty}e^{-\rho t}\bar{q}(t)\bar{v}(t) dt,
  \end{aligned}$$
  where $\check{s}_1=[1,0]\check{s}$ and  $\check{g}(0)=-\int_0^{\infty}e^{-\rho t}\check{s}_1^2dt$. Let $z_s=[z_s^{(1)},
  z_s^{(2)}, z_s^{(3)}, z_s^{(4)}, z_s^{(5)}]^T$ and
  $$\zeta_i=[\zeta_i^{(1)},\zeta_i^{(2)}, \zeta_i^{(3)}, \zeta_i^{(4)}, \zeta_i^{(5)} ]^T,\ i=1,2,3.$$
%  Then
%$$\begin{aligned}
%J_{\rm soc}^{(\infty)} %=&2s(0)q_0-\frac{b^2}{r}\int_0^{\infty}e^{-\rho t}s^2(t)dt
%%=&2q_0[1,0](z^{\ddag}+a_1\zeta_1^{\ddag}+a_2\zeta_2^{\ddag}+a_3\zeta_3^{\ddag})\cr
%=&2q_0\Big(z_s^{(2)}+\sum_{i=1}^3\check{a}_i\zeta_i^{(2)}\Big)\cr
%-&\frac{b^2}{r}\int_0^{\infty}e^{-\rho t}\Big(z_s^{(2)}+\sum_{i=1}^3\check{a}_ie^{\check{\lambda}_it}\zeta_i^{(2)}\Big)dt\cr
%+&\int_0^{\infty}e^{-\rho t}\Big(z_s^{(4)}+\sum_{i=1}^3\check{a}_ie^{\check{\lambda}_it}\zeta_i^{(4)}
%\Big)\Big(z_s^{(5)}+\sum_{i=1}^3\check{a}_ie^{\check{\lambda}_it}\zeta_i^{(5)}\Big)
% dt.
%\end{aligned}$$
From (\ref{eq10d}), we obtain $\check{\lambda}_1=-\alpha$ and $\zeta_1=[0,0,0,0,1]^T$, which implies
$\zeta_1^{(2)}=0,\ \zeta_1^{(4)}=0,\ \zeta_1^{(5)}=1 .$
Thus,
$$\begin{aligned}
 & J_{\rm soc}^{(\infty)} \\
=&2q_0(z_s^{(2)}+\check{a}_2\zeta_2^{(2)}+\check{a}_3\zeta_3^{(2)})-\frac{b^2}{r}\Big[
\frac{(z_s^{(2)})^2}{\rho}+\frac{\check{a}_2^2(\zeta_2^{(2)})^2}{\rho-2\check{\lambda}_2}\\
+&\frac{\check{a}_3^2(\zeta_3^{(2)})^2}{\rho-2\check{\lambda}_3}
+\frac{2\check{a}_2z_s^{(2)}\zeta_2^{(2)}}{\rho-\check{\lambda}_2}+\frac{2\check{a}_3z_s^{(2)}\zeta_3^{(2)}}{\rho-\check{\lambda}_3}
+\frac{2\check{a}_2\check{a}_3\zeta_2^{(2)}\zeta_3^{(2)}}
{\rho-\check{\lambda}_2-\check{\lambda}_3}\Big]\cr
+&\Big[ \frac{z_s^{(4)}z_s^{(5)}}{\rho}+\sum_{i=1}^3\frac{\check{a}_iz_s^{(4)}\zeta_i^{(5)}}{\rho-\check{\lambda}_i}
+\frac{\check{a}_2\zeta_2^{(4)}z_s^{(5)}}{\rho-\check{\lambda}_2}
+\sum_{i=1}^3\frac{\check{a}_2\check{a}_i\zeta_2^{(4)}\zeta_i^{(5)}}{{\rho-\check{\lambda}_2}-\check{\lambda}_i}
\cr
+&\frac{\check{a}_3\zeta_3^{(4)}z_s^{(5)}}{\rho-\check{\lambda}_3}
+\sum_{i=1}^3\frac{\check{a}_3\check{a}_i\zeta_3^{(4)}\zeta_i^{(5)}}{\rho-\check{\lambda}_3-\check{\lambda}_i}
\Big]
%+&\Big[ \frac{z_s^{(4)}z_s^{(5)}}{\rho}+\frac{\check{a}_1z_s^{(4)}\zeta_1^{(5)}}{\rho+\alpha}+\frac{\check{a}_2z_s^{(4)}\zeta_2^{(5)}}{\rho-\check{\lambda}_2}
%+\frac{\check{a}_3z_s^{(4)}\zeta_3^{(5)}}{\rho-\check{\lambda}_3}\cr
%+&\frac{\check{a}_2\zeta_2^{(4)}z_s^{(5)}}{\rho-\check{\lambda}_2}+\frac{\check{a}_1\check{a}_2\zeta_1^{(4)}\zeta_2^{(5)}}{\rho-\check{\lambda}_2+\alpha}
%+\frac{\check{a}_2^2\zeta_2^{(4)}\zeta_2^{(5)}}{\rho-2\check{\lambda}_2}
%+\frac{\check{a}_2\check{a}_3\zeta_2^{(4)}\zeta_3^{(5)}}{\rho-\check{\lambda}_2-\check{\lambda}_3}\cr
%+&\frac{\check{a}_3\zeta_3^{(4)}z_s^{(5)}}{\rho-\check{\lambda}_3}+\frac{\check{a}_1\check{a}_3\zeta_1^{(4)}\zeta_3^{(5)}}{\rho-\check{\lambda}_3+\alpha}
%+\frac{\check{a}_2\check{a}_3\zeta_2^{(4)}\zeta_3^{(5)}}{\rho-\check{\lambda}_2-\check{\lambda}_3}
%+\frac{\check{a}_3^2\zeta_3^{(4)}\zeta_3^{(5)}}{\rho-2\check{\lambda}_3}
%\Big]
.
\end{aligned}$$

The comparison of the two costs is shown in Tab. 1 and Fig. 5.
It can be seen that the asymptotic Nash cost $J_{\rm Nash}^{(\infty)}$ is greater than the optimal social cost $J_{\rm soc}^{(\infty)}$, and {more so when $\alpha$ increases}. %"Collusion results in the smaller output and the higher price"
{This illustrates that the collusion of firms leads to rise in the price and drop in outputs and costs \cite[p. 467]{V93}.}

  \begin{table}
  \begin{center}

%\begin{tabular}{|c|c|c|c|c|c|c|c|c|}\hline
% $\alpha$&0.1& 0.2& 0.5&1&5\\ \hline\hline
% $J_{opt}$&   -0.3107&   -2.0859&   -5.4656&  -7.9532&  -11.1061\\ \hline
%$J_{opt}^{(\infty)}$&      -0.7751&  -3.3162&
%   -8.5858&-12.7198&   -18.0192
%\\ \hline
%\end{tabular}
\begin{tabular}{|c|c|c|c|c|c|c|c|c|}\hline
 $\alpha$&0.1& 0.2& 0.5&1&5\\ \hline\hline
 $ J_{\rm Nash}^{(\infty)}$&   -0.31&   -2.09&   -5.47&  -7.95&  -11.11\\ \hline
$ J_{\rm soc}^{(\infty)}$&      -0.78&  -3.32&
   -8.59&-12.72&   -18.02
\\ \hline
\end{tabular}
{\small Tab. 1: Optimal costs under game and social solutions}
  \end{center}
   \end{table}

 \vskip 0.3cm
\begin{center}
 \hspace{0.8cm}$\epsfig{figure=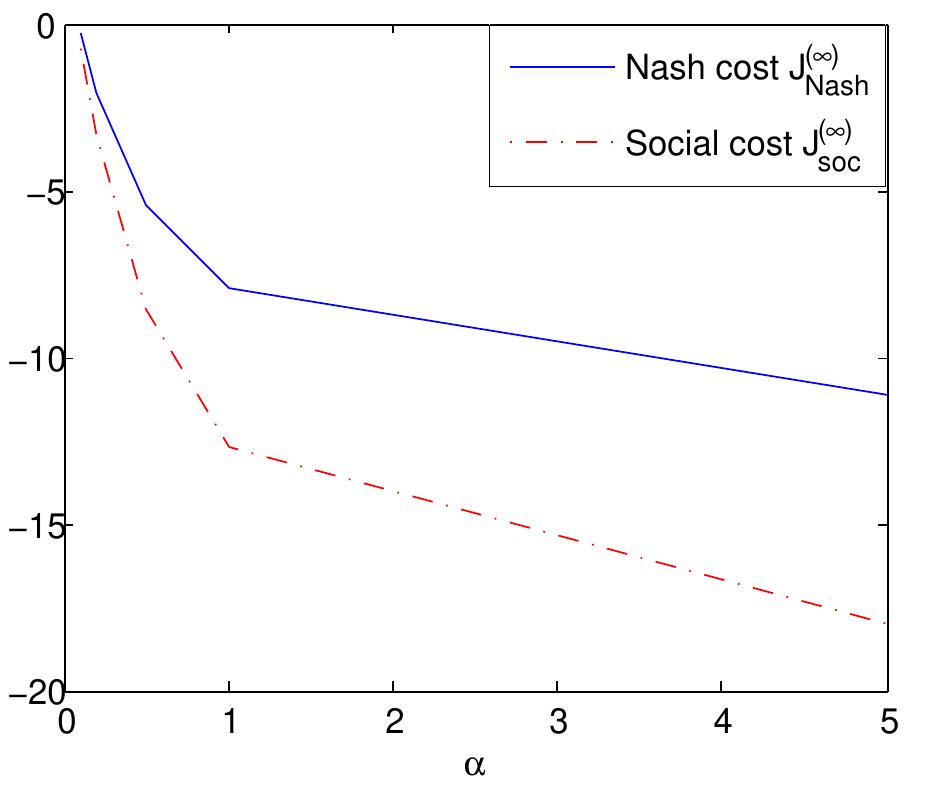,height=5.5cm}$\\
 {\small Fig. 5: Curves of $ J_{\rm Nash}^{(\infty)}$ and $ J_{\rm soc}^{(\infty)}$ with respect to $\alpha$ }
\end{center}

\section{Concluding Remarks}
This paper studies dynamic production output optimization with sticky prices and adjustment costs based on the mean field control methodology. By consistent mean field approximations, we first present the Nash solution for noncooperative firms, and then give the social solution, where agents cooperate to optimize the social cost. The two sets of decentralized strategies are shown to approximate Nash equilibria and social optima, respectively. For future work, it is of interest to consider dynamic production output competition with noisy sticky prices.

\bibliographystyle{plain}

\begin{ack}                               % Place acknowledgements
The first author was supported by the National Natural Science Foundation of China
under Grants 61403233, 61503218. This work was partially conducted while he
was a Research Associate in School of Mathematics and Statistics, Carleton University, Ottawa, Canada during Nov 2014-Apr 2015. The second author was supported by a Discovery Grant of the Natural Sciences and Engineering
Research Council (NSERC) of Canada.
\end{ack}

%\bibliographystyle{plain}        % Include this if you use bibtex
%\bibliography{autosam}           % and a bib file to produce the
                                 % bibliography (preferred). The
                                 % correct style is generated by
                                 % Elsevier at the time of printing.

%\begin{thebibliography}{99}     % Otherwise use the
                                 % thebibliography environment.
                                 % Insert the full references here.
                                 % See a recent issue of Automatica
                                 % for the style.
%  \bibitem[Heritage, 1992]{Heritage:92}
%     (1992) {\it The American Heritage.
%     Dictionary of the American Language.}
%     Houghton Mifflin Company.

\appendix

%\section{A summary of Latin grammar}
\section*{Appendix A: Proof of Theorem \ref{thm3}}

\def\theequation{A.\arabic{equation}}
\setcounter{equation}{0}

\def\thelemma{A.\arabic{lemma}}
\setcounter{lemma}{0}
\def\theproposition{A.\arabic{proposition}}
\setcounter{proposition}{0}

\begin{proposition}\label{prop1}
Under \textbf{A1)}-\textbf{A3)}, $\bar{J}_i(u_i)$ is coercive in $u_i$, i.e., there exists a constant $C$ depending on $p_0$, $q_0$ and a constant $\epsilon_0>0$ such that
\begin{equation}\label{eqa1}
  \bar{J}_i(u_i)\geq \epsilon_0 \|u_i\|_{\rho}^2 -C,
  \end{equation}
where $\|u_i\|_{\rho}^2=E \int_0^{\infty} e^{-\rho t } |u_i(t)|^2dt$.
\end{proposition}

\emph{Proof.}  %By \cite{lim1999stochastic},
%$\bar{J}_i(u_i)$ is strictly convex in $u_i$ if and only if ${J}^{\prime}_i(u_i)$ is positive definite, where
%$$J_i^{\prime}(u_i)=E\int_0^{\infty} e^{-\rho t}r {u}_i^2dt.$$
%It is straightforward to verify $J^{\prime}_i(u_i)$ is positive definite. In fact, we may use a direct method for verifying the definition of convexity to get that  $J^{\prime}_i(u_i)$ is convex.
{From (\ref{eq2}),
\begin{align*}
 %\label{eq7}
 q_i(t)=&q_i(0)e^{-\mu t}
+\int_0^t e^{-\mu (t-\tau)}\big[b_iu_i(\tau)d\tau+{\sigma}dw_i(\tau)\big].
\end{align*}
{Thus, by Cauchy¡¯s inequality and \textbf{A3)},
\begin{align}\label{eq11a}
 &E\int_0^{\infty} e^{-\rho t}|q_i(t)|dt\cr
%\leq& E \int_0^{\infty} e^{-\rho t} |q^{(N)}(0)|e^{-\mu t}dt\cr
%  &+E\int_0^{\infty} e^{-\rho t} \int_0^t e^{-\mu (t-\tau)}\Big[\big|bu^{(N)}(\tau)\big|d\tau+\Big|\frac{\sigma}{N}\sum_{i=1}^Ndw_i(\tau)\Big|\Big]dt \cr
  \leq & C+E\int_0^{\infty} e^{-(\rho+\mu) t} \int_0^t e^{\mu \tau}|b_iu_i(\tau)|d\tau dt\cr
  =& C+E\int_0^{\infty} e^{\mu \tau}|b_iu_i(\tau)| \int_\tau^{\infty} e^{-(\rho+\mu) t} dtd\tau\cr
  =&C+\frac{1}{\rho+\mu}E\int_0^{\infty} e^{-\rho \tau}|b_iu_i(\tau)| d\tau\cr
  \leq &C+\delta_1 E\int_0^{\infty} e^{-\rho \tau}|u_i(\tau)|^2 d\tau,
\end{align}
where $\delta_1$ is a sufficiently small positive number.} Note $\bar{q}\in C_b([0,\infty),\mathbb{R})$.  By (\ref{eq5a}), it follows that $\bar{p}\in C_b([0,\infty),\mathbb{R})$. From this together with (\ref{eq11a}) and (\ref{eq5}), we obtain that there exists a constant $\epsilon_0>0$ such that
(\ref{eqa1}) holds.}
%
%Furthermore, we may get that there exists a small constant $\epsilon_0>0$ such that
%$$\bar{J}_{i,\epsilon_0}(u_i)\stackrel{\Delta}{=}\bar{J}_i(u_i)-2\epsilon_0\|u_i\|^2_{\rho} $$ is convex. By Lemma \ref{lem1a} in \cite{huang2010large}, there exists a constant $C$ dependent of $p_0$ and $q_0$ such that $\sup_{\|u_i\|_{\rho}\leq 1}\big|\bar{J}_{i,\epsilon_0}(u_i)\big|\leq C$. %It can be verified that $\bar{J}_{i,\epsilon_0}(u_i)$ is convex since $\epsilon_0$ is small.
%For $\|u_i\|_{\rho}\geq 1$, it follows from the convexity of  $\bar{J}_{i,\epsilon_0}(u_i)$ that
%$$\begin{aligned}
%\bar{J}_{i,\epsilon_0}\Big(\frac{u_i}{\|u_i\|_{\rho}}\Big)&\leq \frac{1}{\|u_i\|_{\rho}}\bar{J}_{i,\epsilon_0}({u_i})+\frac{\|u_i\|_{\rho}-1}{\|u_i\|_{\rho}}\bar{J}_{i,\epsilon_0}({0})\\
%&\leq
%\frac{1}{\|u_i\|_{\rho}}\bar{J}_{i,\epsilon_0}({u_i})+C,
%\end{aligned}$$
%which implies $\bar{J}_{i,\epsilon_0}({u_i})\geq -2C \|u_i\|_{\rho}.$
%Hence for any $u_i$, $$\bar{J}_{i,\epsilon_0}({u_i})\geq -C(2\|u_i\|_{\rho}+1),$$
%which further leads to
%$$\bar{J}_i(u_i)\geq \epsilon_0 \|u_i\|_{\rho}^2 -C.$$
\rightline{$\Box$}

\begin{lemma}\label{lem1a}
  For any $i=1,\cdots,N$, there exists a constant $C$ such that
  \begin{equation}\label{eq25aa}
    \|q_i\|^2_{\rho}\leq C\|u_i\|^2_{\rho}+C.
      \end{equation}
\end{lemma}
\emph{Proof.} Denote $q_{i,\rho}=e^{-\frac{\rho}{2}}q_i $ and $u_{i,\rho}=e^{-\frac{\rho}{2}}u_i $. By Ito's formula, we have
$$dq_{i,\rho}=-(\mu+\frac{\rho}{2})q_{i,\rho}dt+u_{i,\rho}dt+e^{-\frac{\rho}{2}}dw_i. $$
From this it follows that $$E\int_0^{\infty}|q_{i,\rho}(t)|^2 dt\leq C+E\int_0^{\infty}\left|\int_0^{t}e^{-(\mu+\frac{\rho}{2})(t-\tau)}u_{i,\rho}(\tau)d\tau\right|^2dt.$$
By the argument in the proof of Lemma A.1 in \cite{huang2010large}, we get (\ref{eq25aa}). $\qquad \Box$

\begin{lemma}\label{lem2}
  For the problem (\ref{eq1})-(\ref{eq2}), assume that \textbf{A1)}-\textbf{A4)} hold. For $u_i\in {\mathcal U}_c$, if $J(u_i, \hat{u}_{-i})\leq C_1$, then there exist an integer $N_0$ and a constant $C_2$  such that for all $N\geq N_0$, $\|u_i\|^2_{\rho}\leq  C_2$.
\end{lemma}

\emph{Proof.} %Denote ${\bar{J}}_i(u_i)$ by ${\bar{J}}^*_i(u_i)$ when $\bar{p}$ is replaced by $\bar{{p}}^*$.
For $u_i\in {\mathcal U}_c$, we have
\begin{eqnarray}\label{eq26a}
{J}_i(u_i, \hat{u}_{-i})
&=& {\bar{J}}_i(u_i)- E\int_0^{\infty}e^{-\rho t} [p(t)-\bar{p}(t)] q_i(t)dt\cr
   &\hspace*{-0.2cm}\geq\hspace*{-0.2cm}& {\bar{J}}_i(u_i)-\Big[E\int_0^{\infty}e^{-\rho t} |p(t)-\bar{p}(t)|^2dt\cr
    &&\cdot E\int_0^{\infty}e^{-\rho t}|q_i(t)|^2dt
   \Big]^{1/2}\cr
   &\hspace*{-0.2cm}\stackrel{\Delta}{=}\hspace*{-0.2cm}&{\bar{J}}_i(u_i)-I_1.
  \end{eqnarray}
 % \begin{eqnarray*}
%  \bar{J}_i(u_i)&=& J(u_i)+ E\int_0^{\infty}e^{-\rho t} [p(t)-\bar{p}(t)] q_i(t)dt\cr
%   &\leq& J(u_i)+\Big[E\int_0^{\infty}e^{-\rho t} |p(t)-\bar{p}(t)|^2dt \cdot E\int_0^{\infty}e^{-\rho t}|q_i(t)|^2dt
%   \Big]^{1/2}\cr
%    &\leq&C+o(1).
%  \end{eqnarray*}
By  {Proposition} \ref{prop1}, $\bar{J}_i(u_i)$ is coercive, i.e.,
\begin{equation}\label{eq27a}
  \bar{J}_i(u_i) \geq \epsilon_0 \|u_i\|_{\rho}^2 -C.
  \end{equation}
We now estimate $I_1$. Notice $p(t)-\bar{p}(t)=p(t)-\hat{p}(t)+\hat{p}(t)-\bar{p}(t)$.
We have
\begin{eqnarray*}
d\left[\begin{array}{c}
{p}-\hat{p}\\ {q}^{(N)}-\hat{q}^{(N)}
\end{array}\right]&=&{G}
 \left[\begin{array}{c}
p- \hat{p}\\ {q}^{(N)}-\hat{q}^{(N)}
\end{array}\right]dt +\left[\begin{array}{c}
 0\\
 \frac{b_i}{N}(u_i-\frac{b_is}{r})dt
\end{array}\right]\hspace*{-0.1cm},
\end{eqnarray*}
where $G$ is defined in (\ref{eq22bb}).
%$${G}=\left[\begin{array}{cc}
%-\alpha&-\alpha \\
%0&-\mu
%\end{array}\right].$$
%Then
%\begin{align*}
% & \left[\begin{array}{c}
%{p}-\hat{p}\\ {q}^{(N)}-\hat{q}^{(N)}
%\end{array}\right]
%=\int_0^t e^{{G}( t-\tau)} \left[\begin{array}{c}
% 0\\
%  \frac{b_i}{N}(u_i-\frac{b_is}{r})
%\end{array}\right]d\tau.
%\end{align*}
Noticing that ${G}$ is Hurwitz, by basic estimates as in the proof of in  Lemma \ref{lem1a}, we obtain
\begin{equation}\label{eq27b}
\int_0^{\infty}e^{-\rho t} |p(t)-\hat{p}(t)|^2dt\leq \frac{C}{N^2}+\frac{C}{N^2}\|u_i\|^2_{\rho}.
\end{equation}
From this together with Theorem \ref{thm2}, it follows that
$$E\int_0^{\infty}e^{-\rho t} |p(t)-\bar{p}(t)|^2dt\leq C+\frac{C}{N^2}\|u_i\|^2_{\rho},  $$
which with Lemma \ref{lem1a} implies $$|I_1|\leq C+\frac{C}{N}\|u_i\|^2_{\rho}. $$
Combining this together with (\ref{eq26a}) and (\ref{eq27a}) yields
$$C_1\geq J_i(u_i,\hat{u}_{-i})\geq (\epsilon_0-\frac{C}{N})\|u_i\|^2_{\rho}-2C.$$
Let $N_0=\inf\{m\in \mathbb{Z}|m>C/\epsilon_0\}$.
%[\frac{C}{\epsilon_0}]+1$.
 From this inequality, we obtain that there exists a constant $C_2$ such that for all $N\geq N_0$,
$$\|u_i\|^2_{\rho}\leq \frac{N(C_1+2C)}{N\epsilon_0-C}\leq C_2.$$
\rightline{$\Box$}
%%%%%%%%%%%%%%%%%%%%%%%%%%%%%%%%%%%%%%%%%%%%%%%%%%%%%%%%%%%%%%%%%%%%%%%%%%%%%%%%
%\section*{APPENDIX}
%

\begin{lemma}\label{lem3a}
Under \textbf{A1)}-\textbf{A4)},
$$|\bar{J}_i(\hat{u}_i)-J_i(\hat{u}_i,\hat{u}_{-i})|= O(\varepsilon_N+\frac{1}{\sqrt{N}}).$$
%where $\varepsilon_N=E\int_0^{\infty}e^{-\rho t} |\hat{p}(t)-\bar{p}^*(t)|^2dt=o(1)$.
\end{lemma}

\emph{Proof. } %It follows by Theorem \ref{thm2} that$$  E\int_0^{\infty}e^{-\rho t} |\hat{p}(t)-\bar{p}(t)|^2dt\leq O(\varepsilon_N^2+\frac{1}{N}).$$
    By Schwarz's inequality,
   \begin{align*}
&|\bar{J}_i(\hat{u}_i)-J_i(\hat{u}_i,\hat{u}_{-i})|\\
=&\int_0^{\infty}e^{-\rho t}[\hat{p}(t)\hat{q}_i(t)-\bar{p}(t)\hat{q}_i(t)]dt\\
\leq &\Big[\int_0^{\infty}e^{-\rho t}|\hat{q}_i(t)|^2dt \cdot\int_0^{\infty}e^{-\rho t} |\hat{p}(t)-\bar{p}(t))|^2dt\Big]^{\frac{1}{2}}\\
\leq &O(\varepsilon_N+\frac{1}{\sqrt{N}}),
   \end{align*}
where the last inequality follows from Theorem \ref{thm2}. \hfill $\Box$

\emph{Proof of Theorem \ref{thm3}.} It suffices to show the first inequality in (\ref{eq26n}) under $J_i({u}_i, \hat{u}_{-i})\leq C_1$. By Lemma \ref{lem2}, $J_i({u}_i, \hat{u}_{-i})\leq C_1$
implies
\begin{equation}\label{eq30}
\|u_i\|^2_{\rho}\leq C_2.
\end{equation}
It follows from (\ref{eq26a}) that
\begin{align}\label{eq30b}
{J}_i(u_i, \hat{u}_{-i})={\bar{J}}_i(u_i)-I_1,
  \end{align}
where $$I_1=\Big[E\int_0^{\infty}e^{-\rho t} |p(t)-\bar{p}(t)|^2dt \cdot E\int_0^{\infty}e^{-\rho t}|q_i(t)|^2dt
   \Big]^{1/2}.$$
By (\ref{eq27b}), (\ref{eq30}) and Theorem \ref{thm2},
$$E\int_0^{\infty}e^{-\rho t} |p(t)-\bar{p}(t)|^2dt\leq O(\varepsilon_N^2+\frac{1}{{N}}),$$
which with Lemma \ref{lem1a} further implies
$ |I_1|\leq O(\varepsilon_N+\frac{1}{\sqrt{N}}) .$
Thus, by (\ref{eq30b}) it follows that
\begin{align}
{J}_i(u_i, \hat{u}_{-i})=&\hspace*{0.05cm} {\bar{J}}_i(u_i)-O(\varepsilon_N+\frac{1}{\sqrt{N}})\cr \label{eqA9}
\geq&\hspace*{0.05cm}  {\bar{J}}_i(\hat{u}_i)-O(\varepsilon_N+\frac{1}{\sqrt{N}}).
\end{align}
On the other hand, by Lemma \ref{lem3a} we have
\begin{equation}\label{eqA10}
  {\bar{J}}_i(\hat{u}_i)\geq J_i(\hat{u}_i,\hat{u}_{-i})-O({\varepsilon_N}+\frac{1}{\sqrt{N}}).
\end{equation}
Thus, the first inequality in (\ref{eq26n}) follows from (\ref{eqA9}) and (\ref{eqA10}).
  \hfill $\Box$

\section*{Appendix B: Proof of Theorem \ref{thm4}}

\def\theequation{B.\arabic{equation}}
\setcounter{equation}{0}

\def\thelemma{B.\arabic{lemma}}
\setcounter{lemma}{0}
\def\theproposition{B.\arabic{proposition}}
\setcounter{proposition}{0}
%
%
%To prove Theorem \ref{thm4}, we need some lemmas.

\begin{lemma}\label{lem1}
    Assume that \textbf{A1)}-\textbf{A3)} and \textbf{A5)} hold. For the set of strategies $\{\check{u}_i,i=1,\cdots,N\}$, we have
       \begin{align}\label{eq20}
 &\sup_{t\geq 0, N\geq 1} E\big\{|\check{p}(t)|^2+|\check{q}_i(t)|^2+|\check{q}^{(N)}(t)|^2+|\check{v}_i(t)|^2\big\}\cr
        &\leq C_0,\\\label{eq20b}
&\hspace*{-0.1cm}\sup_{t\geq 0}E \big\{|\check{p}(t)-\bar{p}(t)|^2+|\check{q}^{(N)}(t)-\bar{q}(t)|^2+|\check{v}^{(N)}(t)-\bar{v}(t)|^2\big\}\cr
    &
=O(\varepsilon_N^2+{1}/{N}),
     \end{align}
where $[\bar{q}, \bar{v}]^T=\int_{\Theta}y_{\theta}dF(\theta)$ is given by (\ref{eq9})-(\ref{eq9g}), and $\check{v}^{(N)}=\frac{1}{N}\sum_{i=1}^N \check{v}_i$.

\end{lemma}

\emph{Proof.} It follows from (\ref{eq18}) and (\ref{eq19}) that
\begin{eqnarray*}
\left[\begin{array}{c}
d\check{q}_i\\d\check{v}_i
\end{array}\right]&=&
A
\left[\begin{array}{c}
\check{q}_i\\\check{v}_i
\end{array}\right]dt-\frac{1}{r}B_iB_i^T\check{s}dt+\left[\begin{array}{c}
dw_i\\0
\end{array}\right].
\end{eqnarray*}
 Notice that $A$ is Hurwitz and $\check{s}\in C_b([0,\infty),\mathbb{R}^2).$ By elementary linear estimates and \textbf{A3)} we can show that there exists a constant $C_0$ independent of
 $(i,N)$ such that
$$\sup_{t\geq 0}  E\big\{|\check{q}_i(t)|^2+|\check{v}_i(t)|^2\big\}\leq C_0,$$
which further gives $\sup_{t\geq 0}E\big\{|\check{q}^{(N)}(t)|^2\big\}\leq C_0$.
This together with (\ref{eq17}) leads to (\ref{eq20}).
By (\ref{eq18}) and (\ref{eq19}), we have
\begin{align}\label{eq21}
\left[\begin{array}{c}
d\check{q}^{(N)}\\d\check{v}^{(N)}
\end{array}\right]=&
A
\left[\begin{array}{c}
\check{q}^{(N)}\\
\check{v}^{(N)}
\end{array}\right]dt-\frac{1}{Nr}\sum_{i=1}^NB_iB_i^T\check{s}dt\cr
&
+\left[\begin{array}{c}
\frac{\sigma}{N}\sum_{i=1}^N dw_i\\0
\end{array}\right].
\end{align}
 Denote $\xi=[\check{q}^{(N)}-\bar{q}, \check{v}^{(N)}-\bar{v} ]^T$. It follows from (\ref{eq9}) and (\ref{eq21}) that
\begin{align*}
  d\xi=&A\xi dt+\left[\frac{1}{r}\int_{\Theta} B_{\theta}B_{\theta}^T\check{s}dF(\theta)  -\frac{1}{Nr}\sum_{i=1}^NB_iB_i^T\check{s}\right]dt\cr
&  +\left[\begin{array}{c}
\frac{\sigma}{N}\sum_{i=1}^N dw_i\\0
\end{array}\right],
\end{align*}
where $\xi(0)=[\check{q}^{(N)}(0)-{q}_0, 0]^T$.
%This leads to
%\begin{align}
%  \label{eq22}
%  \xi(t)&=e^{At}\xi(0)+\int_0^t e^{A(t-\tau)}\bar{\Delta}_sd\tau\cr
%  &+\int_0^t e^{A(t-\tau)} \frac{\sigma}{N}\sum_{i=1}^N dw_i(\tau),
%\end{align}
%where
%$$\bar{\Delta}_s=\frac{1}{r}\int_{\Theta} B_{\theta}B_{\theta}^T\check{s}dF(\theta)  -\frac{1}{Nr}\sum_{i=1}^NB_iB_i^T\check{s}. $$
Since $A$ is Hurwitz,
 %and
%$$ \quad \bar{\Delta}_s=O(\varepsilon_N),$$
%$$\begin{aligned}
%  E\|\xi(0)\|^2
% \leq\frac{\max_{1\leq i\leq N}Eq_i^2(0)-q_0^2}{N}
%=O(\frac{1}{N}).
%\end{aligned}$$
%Hence, it follows from (\ref{eq22}) that
\begin{align}\label{eq23}
& E|\check{q}^{(N)}(t)-\bar{q}(t)|^2+E|\check{v}^{(N)}(t)-\bar{v}(t)|^2\cr
=&O(\varepsilon_N^2+{1}/{N}).
     \end{align}
By (\ref{eq9}) and (\ref{eq17}) we have
\begin{align*}
\frac{d(\check{p}-\bar{p})}{dt}&=-\alpha(\check{p}-\bar{p})-\alpha(\check{q}^{(N)}-\bar{q}),\cr
\check{p}(0)&=\bar{p}(0),
  \end{align*}
which leads to
$$\check{p}(t)-\bar{p}(t)=-\alpha\int_0^te^{-\alpha(t-\tau)}[\check{q}^{(N)}(\tau)-\bar{q}(\tau)]d\tau.$$
This together with (\ref{eq23}) gives (\ref{eq20b}).
  \hfill $\Box$

\begin{lemma}\label{lem3}
  There exist constants $C_1>0,C_2>0$ and $C_3>0 $ such that
  $$J^{(N)}_{\rm soc}(u)\geq C_1\int_0^{\infty} e^{-\rho t}p^2dt+\frac{C_2}{N}\int_0^{\infty} e^{-\rho t}\sum_{i=1}^Nu_i^2 dt-C_3.$$
\end{lemma}

\emph{Proof.} {By \textbf{A3)}, assume $|b_i|\leq \hat{b},\ i=1,\cdots,N$. % Denote $u^{(N)}=\frac{1}{N}\sum_{i=1}^Nu_i$.
From (\ref{eq2}), we compute
  %$$dq^{(N)}
%=-\mu q^{(N)}dt+bu^{(N)}dt+\frac{\sigma}{N}\sum_{i=1}^Ndw_i,
%$$
%which yields
\begin{align*}
 %\label{eq7}
 \Big|q^{(N)}(t)\Big|\leq&\Big|q^{(N)}(0)e^{-\mu t}\Big|
+\int_0^t e^{-\mu (t-\tau)}\Big|\frac{1}{N}\sum_{i=1}^Nb_iu_i(\tau)\Big|d\tau\cr
&+\left|\frac{\sigma}{N}\sum_{i=1}^N\int_0^t e^{-\mu (t-\tau)}dw_i(\tau)\right|.
\end{align*}
Thus, by Cauchy's inequality,
\begin{align}\label{eq11}
 &E\int_0^{\infty} e^{-\rho t}|q^{(N)}|dt\cr
%\leq& E \int_0^{\infty} e^{-\rho t} |q^{(N)}(0)|e^{-\mu t}dt\cr
%  &+E\int_0^{\infty} e^{-\rho t} \int_0^t e^{-\mu (t-\tau)}\Big[\big|bu^{(N)}(\tau)\big|d\tau+\Big|\frac{\sigma}{N}\sum_{i=1}^Ndw_i(\tau)\Big|\Big]dt \cr
  \leq & C+E\int_0^{\infty} e^{-(\rho+\mu) t} \int_0^t e^{\mu \tau}\Big|\frac{1}{N}\sum_{i=1}^Nb_iu_i(\tau)\Big|d\tau dt\cr
  =& C+E\int_0^{\infty} e^{\mu \tau}\Big|\frac{1}{N}\sum_{i=1}^Nb_iu_i(\tau)\Big| \int_\tau^{\infty} e^{-(\rho+\mu) t} dtd\tau\cr
  =&C+\frac{1}{\rho+\mu}E\int_0^{\infty} e^{-\rho \tau}\Big|\frac{1}{N}\sum_{i=1}^Nb_iu_i(\tau)\Big| d\tau\cr
  \leq &C+\delta_1 E\int_0^{\infty} e^{-\rho \tau}\Big|\frac{1}{N}\sum_{i=1}^Nb_iu_i(\tau)\Big|^2 d\tau\cr
  \leq &C+\frac{\delta_1 \hat{b}^2}{N} E\int_0^{\infty} e^{-\rho \tau}\sum_{i=1}^Nu_i^2(\tau) d\tau,
\end{align}
where $\delta_1$ is a sufficiently small positive number.} Note $$p(t)=e^{-\alpha t}p(0)+\int_0^te^{-\alpha(t-\tau)}(-\alpha q^{(N)}+\alpha \beta)d\tau.$$
From (\ref{eq11}) we have
\begin{align}\label{eq12}
  &E\int_0^{\infty} e^{-\rho t}|p(t)|dt\cr
\leq& E \int_0^{\infty}e^{-\rho t}\left|e^{-\alpha t}p(0)+\int_0^te^{-\alpha(t-\tau)}(-\alpha q^{(N)}
  +\alpha \beta)d\tau\right|dt\cr
\leq & C+\alpha E\int_0^{\infty}  e^{-(\rho+\alpha)t} \int_0^te^{\alpha \tau}|q^{(N)}(\tau)|d\tau dt\cr
=&C+\frac{\alpha}{\rho+\alpha}E\int_0^{\infty}e^{-\rho \tau}|q^{(N)}(\tau)|d\tau\cr
\leq & C+\frac{\alpha\delta_1 \hat{b}^2}{N(\rho+\alpha)} E\int_0^{\infty} e^{-\rho \tau}\sum_{i=1}^Nu_i^2(\tau) d\tau,
\end{align}
where the fourth line is obtained by an exchange of order of the integration.

Consider the system
\begin{align} \label{eq13}
  \frac{dp}{dt}=&-\alpha p-\alpha \acute{u},\ p(0)=p_0\cr
y=&-p\cr
\acute{u}=&q^{(N)}-\beta.
\end{align}
{By verifying the conditions in the positive real lemma (see e.g., \cite[Lemma 6.2]{Khalil1996}), we obtain that the transfer function of the system (\ref{eq13})
 is positive real, which leads to the passivity of (\ref{eq13}). This implies that there exists a constant $l>0$ such that $\acute{u}y\geq \dot{V}(p)$ where $V(p)\stackrel{\Delta}{=}lp^2$, which further gives}
\begin{align}
 & E\int_0^{\infty} e^{-\rho t}(-p)(q^{(N)}-\beta)dt
 \geq E\int_0^{\infty} e^{-\rho t}d(lp^2)\cr
  %=& l[e^{-\rho t}p^2(t)\left|
%  \begin{array}{c}
%    \infty\\
%    0
%  \end{array}\right. +\rho\int_0^{\infty}p^2e^{-\rho t}dt]
    =& lE\Big[\rho\int_0^{\infty}p^2e^{-\rho t}dt-p^2(0)\Big].
\end{align}
From this together with (\ref{eq11}) and (\ref{eq12}), we have
\begin{align*}
 & J^{(N)}_{\rm soc}(u)\cr
=&   E\int_0^{\infty} e^{-\rho t}\Big[(c-p)q^{(N)}+\frac{r}{N}\sum_{i=1}^Nu_i^2(t)\Big]dt\cr
  =&  E\int_0^{\infty} e^{-\rho t}\Big[(-p)(q^{(N)}-\beta)-\beta p+cq^{(N)}+\frac{r}{N}\sum_{i=1}^Nu_i^2(t)\Big]dt\cr
  \geq & lE\Big[\rho\int_0^{\infty}p^2(t)e^{-\rho t}dt-p^2(0)\Big]\cr
&+\frac{1}{N}(r-\frac{\beta\alpha\delta_1\hat{b}^2}{\rho+\alpha}-c\delta_1 \hat{b}^2)
  E\int_0^{\infty} e^{-\rho t}\sum_{i=1}^Nu_i^2(t) dt-C.
\end{align*}
This completes the proof of the lemma. \hfill $\Box$

{\emph{Proof of Theorem \ref{thm4}}.}
Notice $\inf_{u_i\in {\mathcal U}_c} J^{(N)}_{\rm soc}(u)\leq J^{(N)}(\check{u})\leq C$. It suffices to consider all $u_i\in {\mathcal U}_c$ %$\{u_i, i=1,\cdots,N|u_i\in {\mathcal U}_c\}$
satisfying $$J^{(N)}_{\rm soc}(u)\leq J^{(N)}(\check{u})\leq C. $$
By Lemma \ref{lem3},
\begin{equation}\label{eq24}
  E\int_0^{\infty} e^{-\rho t}\Big(|p|^2+\frac{1}{N}\sum_{i=1}^N|u_i|^2\Big)dt<\infty.
\end{equation}
Let $\tilde{q}_i=q_i-\check{q}_i$, $\tilde{p}=p-\check{p}$, $\tilde{v}_i=v_i-\check{v}_i$ and
 $$\tilde{u}_i={u}_i-\check{u}_i={u}_i+\frac{1}{r}[b_i,0]\check{s}.$$
 By (\ref{eq1}), (\ref{eq2}) and (\ref{eq17})-(\ref{eq19}) we have
\begin{align}
\label{eq25a}\frac{d\tilde{q}_i}{dt}&=-\mu \tilde{q}_i+b_i\tilde{u}_i,\\
\label{eq25} \frac{d\tilde{p}}{dt}&=-\alpha \tilde{p}-\alpha\tilde{q}^{(N)},\\
\label{eq26} \frac{d\tilde{v}_i}{dt}&=-\alpha \tilde{v}_i-\alpha\tilde{q}_i,
   \end{align}
where $\tilde{q}^{(N)}=\frac{1}{N} \sum_{i=1}^N\tilde{q}_i$, and $\tilde{q}_i(0)=\tilde{p}(0)= \tilde{v}_i(0)=0$. Denote $\tilde{v}^{(N)}=\frac{1}{N} \sum_{i=1}^N\tilde{v}_i$.
It follows from (\ref{eq26}) that
$$ \frac{d\tilde{v}^{(N)}}{dt}=-\alpha \tilde{v}^{(N)}-\alpha\tilde{q}^{(N)},\quad  \tilde{v}^{(N)}(0)=0.$$
From this together with (\ref{eq25}) we get $\tilde{v}^{(N)}=\tilde{p}$.
We have
\begin{align}\label{eqb14}
    &J^{(N)}_{\rm soc}(u)\cr
    =&E\int_0^{\infty} e^{-\rho t} \left[(c-p)q^{(N)}+\frac{1}{N}\sum_{i=1}^Nru_i^2\right]dt\cr
    =&J^{(N)}_{\rm soc}(\check{u})+E\int_0^{\infty} e^{-\rho t}\left[\frac{1}{N}\sum_{i=1}^Nr\tilde{u}_i^2-\tilde{p}\tilde{q}^{(N)}\right]dt\cr
    &+ E\int_0^{\infty} e^{-\rho t} \left[(c-\check{p})\tilde{q}^{(N)}-\tilde{p}\check{q}^{(N)} -\frac{1}{N}\sum_{i=1}^N2\tilde{u}_iB_i^Ts\right]dt\cr
   \stackrel{\Delta}{=}&J^{(N)}_{\rm soc}(\check{u})+\tilde{J}^{(N)}_{\rm soc}(\tilde{u})+I^{(N)}.
\end{align}
%where
%\begin{align*}
%\tilde{J}^{(N)}(\tilde{u})&\stackrel{\Delta}{=}\int_0^{\infty} e^{-\rho t}\left[\frac{1}{N}\sum_{i=1}^Nr\tilde{u}_i^2-\tilde{p}\tilde{q}^{(N)}\right]dt\cr
%I^{(N)}&\stackrel{\Delta}{=}\int_0^{\infty} e^{-\rho t} \left[(c-\hat{p})\tilde{q}^{(N)}-\tilde{p}\check{q}^{(N)} -\frac{1}{N}\sum_{i=1}^N2\tilde{u}_iB_i^Ts\right]dt.
%\end{align*}
To complete the proof of the theorem, we will show $J^{(N)}_{\rm soc}(\tilde{u})\geq 0$ and $|I^{(N)}|=O(\frac{1}{\sqrt{N}}+\epsilon_N)$.
For the system
\begin{align*}
\frac{d\tilde{p}}{dt}&=-\alpha \tilde{p}-\alpha u, \quad \tilde{p}(0)=0\cr
y&=-\tilde{p}\cr
u&=\tilde{q}^{(N)},
\end{align*}
by verifying conditions in the positive real lemma (see e.g., \cite{Khalil1996}), we get that the system is passive.
This implies there exists a constant $\tilde{l}>0$ such that
\begin{align}
  E\int_0^{\infty} e^{-\rho t}(-\tilde{p})\tilde{q}^{(N)}dt  \geq \tilde{l}\rho E\Big[\int_0^{\infty}e^{-\rho t}\tilde{p}^2dt\Big]\geq 0,
\end{align}
which with (\ref{eqb14}) leads to $\tilde{J}^{(N)}_{\rm soc}(\tilde{u})\geq 0$.
We have
\begin{equation}\label{eqb15}
\begin{aligned}
I^{(N)}=&\int_0^{\infty} e^{-\rho t} \left[(c-\bar{p})\tilde{q}^{(N)}-\tilde{p}\bar{q} -\frac{1}{N}\sum_{i=1}^N2\tilde{u}_iB_i^Ts\right]dt\cr
+&\int_0^{\infty} e^{-\rho t} \left[( \bar{p}-\check{p})\tilde{q}^{(N)}+\tilde{p}(\bar{q}-\check{q}^{(N)})\right]dt\cr
\stackrel{\Delta}{=}&\zeta_1+\zeta_2.
\end{aligned}
\end{equation}
Applying It\^{o}'s formula to $e^{-\rho t}[\tilde{q}_i, \tilde{v}_i]\check{s}$ and using (\ref{eq9e}) and (\ref{eq9f}), we obtain
\begin{align*}
  &-e^{-\rho T} E\big\{[\tilde{q}_i(T), \tilde{v}_i(T)]\check{s}(T)\big\}
 \cr =&E\int_0^{T} e^{-\rho t} \left[(c-\bar{p})\tilde{q}_i-\tilde{v}_i\bar{q} -2\tilde{u}_iB_i^T\check{s}\right]dt.
  \end{align*}
Noting $\tilde{v}^{(N)}=\tilde{p}$, we get
$$\zeta_1=-\lim_{T\to\infty}\sum_{i=1}^Ne^{-\rho T} E\big\{[\tilde{q}_i(T), \tilde{v}_i(T)]\check{s}(T)\big\}.$$
By (\ref{eq16}) and (\ref{eq24}), it follows that $E|\tilde{u}_i(t)|^2=O(e^{\rho t}).$
This together with (\ref{eq25a}) and (\ref{eq26}) leads to
$E|\tilde{q}_i(t)|^2=O(e^{\rho t})$
and further gives  $E|\tilde{v}_i(t)|^2=O(e^{\rho t}).$
Noticing $\check{s}\in C_b([0,\infty),\mathbb{R}^2),$
 we have $\zeta_1=0$.
As in the proof of Lemma A.1 in \cite{huang2010large}, we use Jenson's inequality %and (\ref{eq7})
to get that
\begin{align}
 & E\int_0^{\infty}e^{-\rho t}|{q}^{(N)}|^2dt\cr
 \leq& C+CE\int_1^{\infty}e^{-\rho t}\Big[\int_0^t e^{-\mu (t-\tau)}\big|\frac{1}{N}\sum_{i=1}^Nb_iu_i(\tau)\big|d\tau\Big]^2dt\cr
 % \leq &C+CE\int_1^{\infty} e^{-\rho t}\int_0^t e^{-\mu (t-\tau)}|u^{(N)}(\tau)|^2\frac{\mu}{1-e^{-\mu t}}d\tau\frac{(1-e^{-\mu t})^2}{\mu^2} dt\cr
  \leq &C+\frac{C}{\mu N(\rho+\mu)}E\int_1^{\infty} e^{-\rho t}\sum_{i=1}^N|u_i(t)|^2dt.
\end{align}
This together with (\ref{eq24}) and (\ref{eq20}) implies
$$ E\int_0^{\infty}e^{-\rho t}|\tilde{q}^{(N)}|^2dt<\infty.$$
By Schwarz's inequality and Lemma \ref{lem1}, we have
\begin{equation}\label{eqb17}
  \int_0^{\infty} e^{-\rho t}( \bar{p}-\check{p})\tilde{q}^{(N)}dt=O(\frac{1}{\sqrt{N}}+\epsilon_N).
  \end{equation}
Furthermore, it follows from Schwarz's inequality, (\ref{eq24}) and Lemma \ref{lem1} that
$$\int_0^{\infty} e^{-\rho t} \tilde{p}(\bar{q}-\check{q}^{(N)})dt=O(\frac{1}{\sqrt{N}}+\epsilon_N).$$
Thus, by (\ref{eqb15}) and (\ref{eqb17}) we have $$|I^{(N)}|=|\zeta_2|=O(\frac{1}{\sqrt{N}}+\epsilon_N).$$
This completes the proof.  \hfill $\Box$   % Each appendix must have a short title.
      % Sections and subsections are supported
                                        % in the appendices.
\end{document}